\documentstyle[amscd]{amsart}
\numberwithin{equation}{section}
\theoremstyle{plain}
\newtheorem{thm}{Theorem}[section]
\newtheorem{cor}[thm]{Corollary}
\newtheorem{lem}[thm]{Lemma}
\newtheorem{prop}[thm]{Proposition}
\begin{document}
\title{$C^*$-algebras associated with 
Hilbert $C^*$-quad modules 
of finite type
}
\author{Kengo Matsumoto}
\address{ 
Department of Mathematics, 
Joetsu University of Education,
Joetsu 943-8512, Japan}
\email{kengo@@juen.ac.jp}
\maketitle
\begin{abstract}
A Hilbert $C^*$-quad module of finite type  
has a multi structure of Hilbert $C^*$-bimodules with two finite bases.
We will construct a $C^*$-algebra from a Hilbert $C^*$-quad module 
of finite type
and prove its universality subject to certain relations among generators.
Some examples of the $C^*$-algebras from Hilbert $C^*$-quad modules of finite type  will be presented. 
\end{abstract}


\def\Zp{{ {\Bbb Z}_+ }}
\def\CTDS{{ ({\cal A}, \rho, \eta, \Sigma^\rho, \Sigma^\eta, \kappa)}}
\def\FKL{{{\cal F}_k^{l}}}
\def\FHK{{{\cal F}_{{\cal H}_\kappa}}}
\def\FU{{{\cal F}_{{\cal H}_\kappa}^{uni}}}
\def\FHLI{{{\cal F}_{[{\frak L}]} }}
\def\FRHO{{ {\cal F}_\rho}}
\def\DRHO{{{{\cal D}_{\rho}}}}
\def\ORHO{{{\cal O}_\rho}}
\def\OETA{{{\cal O}_\eta}}
\def\ORE{{{\cal O}^{\kappa}_{\rho,\eta}}}
\def\OREK{{{\cal O}_{(\rho,\eta;\kappa)}}}
\def\OHK{{{\cal O}_{{\cal H}_\kappa}}}
\def\OH{{{\cal O}_{{\cal H}}}}
\def\FH{{{\cal F}_{{\cal H}}}}
\def\OFH{{{\cal O}_{F({\cal H})}}}
\def\TFH{{{\cal T}_{F({\cal H})}}}
\def\PH{{{\cal P}_{{\cal H}}}}
\def\OU{{{\cal O}_{{\cal H}_\kappa}^{uni}}}
\def\OAK{{{\cal O}_{\A,\kappa}}}
\def\OHRE{{{\widehat{\cal O}}^{\kappa}_{\rho,\eta}}}
\def\FRE{{{\cal F}_{\rho,\eta}}}
\def\DRE{{{\cal D}_{\rho,\eta}}}
\def\M{{{\cal M}}}
\def\N{{{\cal N}}}
\def\H{{{\cal H}}}
\def\HK{{{\cal H}_\kappa}}
\def\EK{{E_\kappa}}
\def\SK{{\Sigma_\kappa}}
\def\S{{{\tt S}}}
\def\T{{{\tt T}}}
\def\K{{{\cal K}}}
\def\A{{{\cal A}}}
\def\B{{{\cal B}}}
\def\BE{{{\cal B}_\eta}}
\def\BR{{{\cal B}_\rho}}
\def\BK{{{\cal B}_\kappa}}
\def\BU{{{\cal B}_\kappa^{uni}}}
\def\Aut{{{\operatorname{Aut}}}}
\def\End{{{\operatorname{End}}}}
\def\Hom{{{\operatorname{Hom}}}}
\def\Ker{{{\operatorname{Ker}}}}
\def\ker{{{\operatorname{ker}}}}
\def\Coker{{{\operatorname{Coker}}}}
\def\id{{{\operatorname{id}}}}
\def\dim{{{\operatorname{dim}}}}
\def\min{{{\operatorname{min}}}}
\def\Max{{{\operatorname{Max}}}}
\def\exp{{{\operatorname{exp}}}}
\def\supp{{{\operatorname{supp}}}}
\def\Proj{{{\operatorname{Proj}}}}
\def\Im{{{\operatorname{Im}}}}
\def\span{{{span}}}
\def\OA{{{\cal O}_A}}
\def\OB{{{\cal O}_B}}

\section{Introduction}

G. Robertson--T. Steger \cite{RoSte}
have initiated a certain study of higher dimensional analogue of Cuntz--Krieger algebras
from the view point of tiling systems of 2-dimensional plane.
After their work, 
A. Kumjian--D. Pask \cite{KP}
have generalized their construction to introduce the notion of 
higher rank graphs and its $C^*$-algebras.
Since then, there have been 
many studies  on these $C^*$-algebras by many authors
(see for example \cite{Dea}, \cite{EGS}, \cite{KP}, \cite{PCHR}, \cite{PRW1}, \cite{RoSte}, etc.).

 In \cite{MaCrelle}, the author has introduced a notion of 
 $C^*$-symbolic dynamical system, 
 which is a generalization of a finite labeled graph,
  a $\lambda$-graph system and an automorphism of a unital $C^*$-algebra.
It is denoted by $(\A, \rho,\Sigma)$
and consists of a finite family $\{ \rho_{\alpha} \}_{\alpha \in \Sigma}$ 
of endomorphisms of a unital $C^*$-algebra 
${\cal A}$ 
such that 
$\rho_\alpha(Z_\A) \subset Z_\A, \alpha \in \Sigma$
and
$\sum_{\alpha\in \Sigma}\rho_\alpha(1) \ge 1$
where
$Z_\A$ denotes the center of $\A$.
It provides a subshift 
$\Lambda_\rho$  
over $\Sigma$
and a Hilbert $C^*$-bimodule 
${\cal H}_{\cal A}^{\rho}$
over $\A$
which gives rise to a  $C^*$-algebra
${\cal O}_\rho$ as a Cuntz-Pimsner algebra 
(\cite{MaCrelle}, cf. \cite{KPW}, \cite{MuSo}, \cite{Pim}).
In \cite{MaPre20111} and \cite{MaPre20112}, 
the author has extended the notion of
$C^*$-symbolic dynamical system to 
$C^*$-textile  dynamical system
which is a higher dimensional analogue of $C^*$-symbolic dynamical system.
A $C^*$-textile  dynamical system
$({\cal A}, \rho, \eta, \Sigma^\rho, \Sigma^\eta, \kappa)$
 consists of two 
$C^*$-symbolic dynamical systems 
$({\cal A}, \rho,  {\Sigma^\rho})$
and 
$({\cal A}, \eta, {\Sigma^\eta})$
with common unital $C^*$-algebra $\A$ and 
commutation relations $\kappa$
between the endomorphisms $\rho_\alpha, \alpha \in \Sigma^\rho$ 
and $\eta_a, a \in \Sigma^\eta$.
A $C^*$-textile dynamical system
 provides  a two-dimensional subshift and 
a multi structure of Hilbert $C^*$-bimodules
that has  multi right actions and 
multi left actions and multi inner products.
Such  a multi structure of Hilbert $C^*$-bimodule is called 
a Hilbert $C^*$-quad module.
In \cite{MaPre20112},
the author has introduced 
a $C^*$-algebra 
associated with the Hilbert $C^*$-quad module
of $C^*$-textile dynamical system.
It is generated by the quotient images 
of creation operators on two-dimensional analogue of Fock Hilbert 
module by module maps of compact operators. 
As a result,
the $C^*$-algebra
has been proved to have a universal property 
subject to certain operator relations 
of generators encoded by structure of 
$C^*$-textile dynamical system (\cite{MaPre20112}).

In this paper, we will generalize the construction of 
the $C^*$-algebras of Hilbert $C^*$-quad modules
of $C^*$-textile dynamical systems.
Let $\A, \B_1, \B_2$ 
be unital $C^*$-algebras.
Assume that $\A$ has unital embeddings into both
$\B_1$ and $\B_2$.
A Hilbert $C^*$-quad module $\H$
over $(\A;\B_1,\B_2)$
is a Hilbert $C^*$-bimodule over $\A$ with 
$\A$-valued right inner product 
$\langle \cdot \mid \cdot \rangle_\A$
which has a multi structure 
of Hilbert $C^*$-bimodules over $\B_i$ with right actions
$\varphi_i$ of $\B_i$ and left actions $\phi_i$ of $\B_i$
and $\B_i$-valued inner products 
$\langle \cdot \mid \cdot \rangle_{\B_i}$
for $i=1,2$
satisfying certain compatibility conditions.
A Hilbert $C^*$-quad module $\H$ is said to be of finite type
if there exist
a finite basis $\{ u_1,\dots,u_M\}$ of $\H$ 
as a Hilbert $C^*$-right module over $\B_1$
and
a finite basis $\{ v_1,\dots,v_N\}$ of $\H$ 
as a Hilbert $C^*$-right module over $\B_2$
 such that
\begin{equation*} 
\sum_{i=1}^M \langle u_i \mid 
\phi_2(\langle \xi \mid \eta \rangle_{\B_2}) u_i \rangle_{\B_1}
=
\sum_{k=1}^N \langle v_k \mid 
\phi_1(\langle \xi \mid \eta \rangle_{\B_1}) v_k \rangle_{\B_2}
=\langle \xi \mid \eta \rangle_\A
\end{equation*}
for $\xi,\eta \in \H$ 
(see \cite{Wat} for the original definition of finite basis
of Hilbert module). 
For a Hilbert $C^*$-quad module,
we will construct a Fock space $F(\H)$ from $\H$,
which is a 2-dimensional analogue to the ordinary Fock space of Hilbert 
$C^*$-bimodules (cf. \cite{Kat}, \cite{Pim}).
We will then define two kinds of creation operators
$s_\xi, t_\xi$ for $\xi \in \H$ on $F(\H)$.  
The $C^*$-algebra on $F(\H)$ generated by them
is denoted by $\TFH$ and called the Toeplitz quad module algebra.
We then define the $C^*$-algebra 
$\OFH$ associated with the Hilbert $C^*$-quad module $\H$ 
by the quotient $C^*$-algebra of $\TFH$ by the ideal generated 
by the finite rank operators.
We will then prove that the $C^*$-algebra 
$\OFH$ for a $C^*$-quad module
$\H$  of finite type has a universal property
in the following way.
\begin{thm}[Theorem \ref{thm:main2}]
Let $\H$ be a Hilbert $C^*$-quad module over
$(\A;\B_1,\B_2)$
of finite type
with
a finite basis $\{ u_1,\dots,u_M\}$ of $\H$ 
as a Hilbert $C^*$-right module over $\B_1$
and
a finite basis $\{ v_1,\dots,v_N\}$ of $\H$ 
as a Hilbert $C^*$-right module over $\B_2$.
Then
the  $C^*$-algebra 
$\OFH$ 
generated by the quotients
$[s_\xi], [t_\xi]$ of the  creation operators 
$s_\xi, t_\xi$ for $\xi \in \H$
on the Fock spaces $F(\H)$
is canonically isomorphic to the universal $C^*$-algebra
$\OH$
generated by operators
$S_1,\dots, S_M, T_1, \dots, T_N$
and elements
$z\in \B_1, w\in \B_2 $
subject to the relations: 
\begin{align*}
\sum_{i=1}^M S_i S_i^* +
 \sum_{k=1}^N T_k T_k^* & =1, 
 \qquad
S_j^* T_l =0,
 \\ 
S_i^* S_j  = \langle u_i \mid u_j \rangle_{\B_1}, 
& \qquad
T_k^* T_l  = \langle v_k \mid v_l \rangle_{\B_2}, \\
z S_j = \sum_{i=1}^M S_i  \langle u_i \mid \phi_1(z) u_j \rangle_{\B_1}, 
& \qquad
z T_l = \sum_{k=1}^N T_k  \langle v_k \mid \phi_1(z) v_l \rangle_{\B_2}, \\
w S_j = \sum_{i=1}^M S_i  \langle u_i \mid \phi_2(w) u_j \rangle_{\B_1},
& \qquad
w T_l = \sum_{k=1}^N T_k  \langle v_k \mid \phi_2(w) v_l \rangle_{\B_2}
\end{align*}
for
$
z \in \B_1, w \in \B_2, i,j=1,\dots, M,
k,l=1,\dots,N$.
\end{thm}
The eight relations of the operators above are called the relations $(\H)$.  
As a corollary we have 
\begin{cor}[Corollary \ref{cor:cor2}]
For a $C^*$-quad module $\H$
of finite type,
the universal $C^*$-algebra
$\OH$
generated by operators
$S_1,\dots, S_M, T_1, \dots, T_N$
and  elements
$z\in \B_1, w\in \B_2 $
subject to the relations
$(\H)$
does not depend on the choice of the finite bases
$\{u_1,\dots,u_M\}$ and
$\{v_1,\dots,v_N\}$.
\end{cor}

The paper is organized in the following way.
In Section 2, we will define Hilbert $C^*$-quad module
and present some basic properties.
In Section 3, we will define a $C^*$-algebra $\OFH$ 
from  Hilbert $C^*$-quad module $\H$ of general type
by using creation operators on Fock Hilbert $C^*$-quad module.
In Section 4, we will study algebraic structure  
of the $C^*$-algebra $\OFH$ 
for a Hilbert $C^*$-quad module $\H$ of finite type.
In Section 5, we will prove, as a main result of the paper,
 that the $C^*$-algebra $\OFH$ 
has the universal property stated as in Theorem 1.1.
A sterategy to prove Theorem 1.1 is to show  
that the $C^*$-algebra $\OFH$ is regarded as a Cuntz-Pimsner algebra 
for a Hilbert $C^*$-bimodule over the $C^*$-algebra generated by 
$\phi_1(\B_1)$ and 
$\phi_2(\B_2)$.
We will then prove the gauge invariant universarity of the $C^*$-algebra
(Theorem \ref{thm:main1}). 
In Section 6, we will present K-theory formulae for the $C^*$-algebra
$\OH$.
In Section 7, we will give examples.
In Section 8, we will formulate higher dimensional anlogue of our situations
and state a generalized proposition of Theorem 1.1 without proof.

Throughout the paper,
we will denote by $\Zp$ 
the set of nonnegative integers 
and
by ${\Bbb N}$
the set of positive integers.

\section{Hilbert $C^*$-quad modules }
Throughout the paper 
we fix three unital $C^*$-algebras 
$\A, \B_1, \B_2$
such that
$\A \subset \B_1, \A \subset \B_2$
with common units.
We assume that there exists a right action $\psi_i$
of $\A$ on $\B_i$ so that
\begin{equation*}
b_i \psi_i(a) \in \B_i \qquad
\text{ for } \quad 
b_i \in \B_i, \, a \in \A, \, i=1,2,
\end{equation*} 
which satisfies
\begin{equation*}
\| b_i \psi_i(a) \| \le \| b_i \| \| a \|, \qquad
b_i \psi_i(a a') = b_i \psi_i(a)\psi_i(a')
\end{equation*}
for
$b_i \in \B_i, \, a, a' \in \A, \, i=1,2.$
Hence $\B_i$ is a right $\A$-module through $\psi_i$ for $i=1,2$.
Suppose that 
$\H$ is a Hilbert $C^*$-bimodule over $\A$,
which has a right action of $\A$, an $\A$-valued 
right inner product $\langle \cdot \mid \cdot \rangle_\A$
and a $*$-homomorphism 
$\phi_\A$ from $\A$
to the algebra of all bounded
adjointable right $\A$-module maps ${\cal L}_\A(\H)$
satisfying
\begin{enumerate}
\renewcommand{\labelenumi}{(\roman{enumi})}
\item $\langle \cdot \mid \cdot \rangle_\A$ is linear in the second variable.

\item $\langle \xi \mid \eta a \rangle_\A = \langle \xi \mid \eta  \rangle_\A a$
for  $\xi, \eta \in \H, a \in \A$.

\item $\langle \xi \mid \eta  \rangle_\A^* 
= \langle \eta \mid \xi  \rangle_\A$
for  $\xi, \eta \in \H$.

\item $\langle \xi \mid \xi \rangle_\A \ge 0$,
and $\langle \xi \mid \xi \rangle_\A = 0$
if and only if $\xi =0$.
\end{enumerate}
A Hilbert $C^*$-bimodule $\H$ over $\A$
is called a Hilbert $C^*$-{\it quad module over} 
$(\A;\B_1,\B_2)$
if  $\H$ has a further structure of a Hilbert $C^*$-bimodule over $\B_i$
for each $i=1,2$
with right action $\varphi_i$ of $\B_i$ and left action $\phi_i$ of $\B_i$
and $\B_i$-valued right inner product 
$\langle \cdot \mid \cdot \rangle_{\B_i}
$
such that 
for $z \in \B_1, w \in \B_2$
$$ \phi_1(z) , \phi_2(w) \in {\cal L}_\A(\H)
$$ 
and
\begin{align}
[\phi_1(z) \xi]\varphi_2(w) = \phi_1(z)[ \xi \varphi_2(w) ], 
& \quad 
[\phi_2(w) \xi]\varphi_1(z) = \phi_2(w)[ \xi\varphi_1(z) ], \\
 \xi \varphi_1(z \psi_1(a)) = [\xi \varphi_1(z)]a, 
& \quad 
\xi \varphi_2(w \psi_2(a)) = [\xi \varphi_2(w)]a, 
\end{align}
for  $\xi \in \H, z\in \B_1, w \in \B_2, a \in \A$
and
\begin{equation}
\phi_\A(a) = \phi_1(a) = \phi_2(a) \qquad 
\text{ for } a \in \A 
\end{equation}
where 
$\A$ is regarded as a subalgebra of $\B_i$.
The left action $\phi_i$ of $\B_i$ 
on $\H$ means that
$\phi_i(b_i)$ for $b_i \in \B_i$
is a bounded adjointable operator with respect to
the inner product $\langle \cdot \mid \cdot \rangle_{\B_i}$
for each $i=1,2$.
The operator 
$\phi_i(b_i)$ for $b_i \in \B_i$
is also adjointable with respect to 
the inner product $\langle \cdot \mid \cdot \rangle_{\A}$.
We assume that 
the adjoint of $\phi_i(b_i)$
with respect to
the inner product $\langle \cdot \mid \cdot \rangle_{\B_i}$
coincides with 
the adjoint of $\phi_i(b_i)$
with respect to
the inner product $\langle \cdot \mid \cdot \rangle_{\A}$.
Both of them coincide with 
$\phi_i(b_i^*)$.
We assume that the left actions $\phi_i$
of $\B_i$ on $\H$ for $i=1,2$ are faithful.
We require the following compatibility conditions between 
the right $\A$-module structure of $\H$ and 
the right $\A$-module structure of $\B_i$
through $\psi_i$:
\begin{equation}
\langle \xi \mid \eta a \rangle_{\B_i}
=\langle \xi \mid \eta  \rangle_{\B_i} \psi_i(a)
\qquad
\text{ for }
\xi, \eta \in \H, a \in \A, \, i=1,2. 
\end{equation}

We further 
assume that $\H$ is a full Hilbert $C^*$-bimodule with respect to 
the three  inner products
$
\langle \cdot \mid \cdot \rangle_\A,
\langle \cdot \mid \cdot \rangle_{\B_1},
\langle \cdot \mid \cdot \rangle_{\B_2}
$
for each.
This means that the $C^*$-algebras generated by elements
$
\{ \langle \xi \mid \eta \rangle_\A \mid \xi,\eta \in \H \},
$
$
\{ \langle \xi \mid \eta \rangle_{\B_1} \mid \xi,\eta \in \H \}
$
and
$
\{ \langle \xi \mid \eta \rangle_{\B_2} \mid \xi,\eta \in \H \}
$
coincide with
$\A$, $\B_1$ and $\B_2$ respectively.

For a vector $\xi \in \H$, denote by
$\| \xi \|_\A, \| \xi \|_{\B_1}, \| \xi \|_{\B_2}$
the norms  
$
\| \langle \xi \mid \xi \rangle_\A \|^{\frac{1}{2}}, 
\| \langle \xi \mid \xi \rangle_{\B_1} \|^{\frac{1}{2}}, 
\| \langle \xi \mid \xi \rangle_{\B_2} \|^{\frac{1}{2}} 
$
induced by the right inner products respectively.
By definition,
$\H$ is complete under the above three norms for each.

\noindent
{\bf Definition.}

(i) 
A Hilbert $C^*$-quad module $\H$ over $(\A;\B_1,\B_2)$
is said to be {\it of general type}
if there exists a faithful completely positive map
$\lambda_i: \B_i \longrightarrow \A$ for $ i=1,2$ 
such that
\begin{align}
\lambda_i(b_i \psi_i(a))&  = \lambda_i(b_i) a
\qquad \text{ for } b_i \in \B_i, a \in \A,  \label{eqn:lambdabpsi}\\ 
\lambda_i(\langle \xi \mid \eta \rangle_{\B_i}) 
& = 
\langle \xi \mid \eta \rangle_{\A}, \qquad
\text{ for }
\xi, \eta \in \H, \, i=1,2. \label{eqn:lambda}
\end{align}

(ii)
A Hilbert $C^*$-quad module $\H$ over $(\A;\B_1,\B_2)$
is said to be {\it of finite type}
if there exist a finite basis 
$\{ u_1,\dots,u_M\}$ of $\H$ as a right Hilbert $\B_1$-module
and a finite basis 
$\{ v_1,\dots,v_N\}$ of $\H$ as a right Hilbert $\B_2$-module,
that is,
\begin{equation}
\sum_{i=1}^M u_i \varphi_1(\langle u_i \mid \xi \rangle_{\B_1})
=
\sum_{k=1}^N v_k \varphi_2(\langle v_k \mid \xi \rangle_{\B_2})
= \xi,
\qquad \xi \in \H \label{eqn:basis}
\end{equation}
such that
\begin{align}
\langle u_i \mid \phi_2(w) u_j \rangle_{\B_1} 
& \in \A, \qquad i,j = 1,\dots,M, \label{eqn:uuA} \\ 
\langle v_k \mid \phi_1(z) v_l \rangle_{\B_2} 
& \in \A, \qquad k,l = 1,\dots,N \label{eqn:vvA}\\
\intertext{for $w \in \B_2$, $z \in \B_1$  and}
\sum_{i=1}^M \langle u_i \mid 
\phi_2(\langle \xi \mid \eta \rangle_{\B_2}) u_i \rangle_{\B_1}
& = \langle \xi \mid \eta \rangle_\A, \label{eqn:uxieta}\\
\sum_{k=1}^N \langle v_k \mid 
\phi_1(\langle \xi \mid \eta \rangle_{\B_1}) v_k \rangle_{\B_2}
& = \langle \xi \mid \eta \rangle_\A \label{eqn:vxieta}
\end{align}
for all $\xi,\eta \in \H$.
Following \cite{Wat}, 
$\{ u_1,\dots,u_M\}$
and  
$\{ v_1,\dots,v_N\}$ 
are called  finite bases of $\H$ respectively.

(iii)
A Hilbert $C^*$-quad module $\H$ over $(\A;\B_1,\B_2)$
is said to be {\it of strongly finite type}
if it is of finite type and there exist a finite basis 
$\{ e_1,\dots,e_{M'} \}$ of $\B_1$ as a right $\A$-module
through $\psi_1\circ \lambda$ 
and a finite basis 
$\{ f_1,\dots,f_{N'} \}$ of $\B_2$ as a right $\A$-module
through $\psi_2\circ \lambda_2$.
This means that the following equalities hold:
\begin{align}
z & 
=\sum_{j=1}^{M'} e_j \psi_1(\lambda_1(e_j^* z)), \qquad z \in \B_1, \\
w & 
=\sum_{l=1}^{N'} f_l \psi_2(\lambda_2(f_l^* w)), \qquad w \in \B_2.
\end{align}

We note that for a Hilbert $C^*$-quad module of general type,
the conditions \eqref{eqn:lambda} imply
$$
\| \langle \xi \mid \xi \rangle_\A \| \le  
\| \lambda_i(1) \|
\| \langle \xi \mid \xi \rangle_{\B_i} \|,
\qquad 
\xi \in \H.
$$ 
Put
$C_i = \| \lambda_i(1) \|^{\frac{1}{2}} >0$
so that
$\| \xi \|_\A \le C_i \| \xi \|_{\B_i}$.
Hence the identity operators
from the Banach spaces
$(\H, \| \cdot \|_{\B_i})$ to 
$(\H, \| \cdot \|_{\A})$
are bounded linear maps.
By the inverse mapping theorem,
there exist constants $C'_i$ such that
$\| \xi \|_{\B_i} \le C'_i \| \xi \|_\A$
for $\xi \in \H$.
  Therefore the three norms 
$\| \cdot \|_\A, \| \cdot \|_{\B_i}, i=1,2,$
induced by the three inner products
$
 \langle \cdot \mid \cdot \rangle_\A,   
 \langle \cdot \mid \cdot \rangle_{\B_i}, i=1,2
$
on $\H$ 
are equivalent to each other.

\begin{lem}
Let $\H$ be a Hilbert $C^*$-quad module $\H$ over $(\A;\B_1,\B_2)$.
If $\H$ is of finite type,
then it is of general type.
\end{lem}
\begin{pf}
Suppose that 
$\H$ is of finite type with finite bases
$\{ u_1,\dots,u_M\}$ of $\H$ as a right Hilbert $\B_1$-module
and 
$\{ v_1,\dots,v_N\}$ of $\H$ as a right Hilbert $\B_2$-module
as above.
We put
\begin{align*}
\lambda_1(z) & =\sum_{k=1}^N \langle v_k \mid \phi_1(z) v_k \rangle_{\B_2} 
\in \A, \qquad z \in \B_1,\\
\lambda_2(w) & =\sum_{i=1}^M \langle u_i \mid \phi_2(w) u_i \rangle_{\B_1} 
\in \A, \qquad w \in \B_2.  
\end{align*}
They give rise to
faithful
completely positive maps 
$\lambda_i: \B_i \longrightarrow \A, i=1,2$.
The equalities
\eqref{eqn:uxieta}
\eqref{eqn:vxieta}
imply that 
\begin{equation}
\lambda_i(\langle \xi \mid \eta \rangle_{\B_i}) 
 = 
\langle \xi \mid \eta \rangle_{\A}, \qquad
\text{ for }
\xi, \eta \in \H, \, i=1,2.
\end{equation}
It then follows that 
\begin{equation*}
\lambda_i(\langle \xi \mid \eta \rangle_{\B_i} \psi_i(a))
 = \lambda_i(\langle \xi \mid \eta a \rangle_{\B_i}) 
 = \langle \xi \mid \eta \rangle_{\A} a 
 = \lambda_i(\langle \xi \mid \eta \rangle_{\B_i} )a,
 \qquad a \in \A.
\end{equation*}
Since $\H$ is full,
the equalities \eqref{eqn:lambdabpsi}
hold.
\end{pf}

\begin{lem}
Suppose that a Hilbert $C^*$-quad module $\H$ of finite type is 
of strongly finite type with 
a finite basis 
$\{ e_1,\dots,e_{M'} \}$ of $\B_1$ as a right $\A$-module
through $\psi_1\circ \lambda_1$ 
and a finite basis 
$\{ f_1,\dots,f_{N'} \}$ of $\B_2$ as a right $\A$-module
through $\psi_2\circ \lambda_2$.
Let
$\{ u_1,\dots,u_M\}$
and  
$\{ v_1,\dots,v_N\}$ 
be finite bases of $\H$
satisfying \eqref{eqn:basis}.
Then two  families 
$
\{ u_i \varphi_1(e_j) \mid i=1,\dots, M,\,  j=1,\dots, M'\} 
$
and
$
\{ v_k \varphi_2(f_l) \mid k=1,\dots, N,\,  l=1,\dots, N'\} 
$
of $\H$ form bases of $\H$ as right $\A$-modules respectively.
\end{lem}
\begin{pf}
For $\xi \in \H$, by the equalities
\begin{equation*}
\xi 
= \sum_{i=1}^M u_i \varphi_1(\langle u_i \mid \xi \rangle_{\B_1}),
\qquad
\langle u_i \mid \xi \rangle_{\B_1} 
= \sum_{j=1}^{M} e_j \psi_1(\lambda_1(e_j^* \langle u_i \mid \xi \rangle_{\B_1})), 
\end{equation*}
it follows that
\begin{align*}
\xi
& =\sum_{i=1}^M u_i \varphi_1(
\sum_{j=1}^{M'} e_j \psi_1(\lambda_1(e_j^* 
\langle u_i \mid \xi \rangle_{\B_1}
))) \\
& =\sum_{i=1}^M \sum_{j=1}^{M'}
u_i \varphi_1(e_j )
 \cdot \lambda_1(e_j^* 
\langle u_i \mid \xi \rangle_{\B_1}) \\
& =\sum_{i=1}^M \sum_{j=1}^{M'}
u_i \varphi_1(e_j )
 \cdot \lambda_1( 
\langle u_i \varphi_1(e_j) \mid \xi \rangle_{\B_1}) \\
& =\sum_{i=1}^M \sum_{j=1}^{M'}
u_i \varphi_1(e_j )
 \cdot  
\langle u_i \varphi_1(e_j) \mid \xi \rangle_{\A}.
\end{align*}
We similarly have
\begin{equation*}
\xi
 =\sum_{k=1}^N \sum_{l=1}^{N'}
v_k \varphi_2(f_l )
 \cdot  
\langle v_k \varphi_2(f_l) \mid \xi \rangle_{\A}.
\end{equation*}
\end{pf}

We present some  examples.

{\bf Examples. }

{\bf 1.}
Let 
$\alpha$,
$\beta$
be automorphisms of a unital 
$C^*$-algebra $\A$
satisfying 
$
\alpha \circ \beta
=
\beta \circ \alpha.
$
Define right actions 
$\psi_i$ of $\A$ on $\B_i$
by
\begin{equation*}
b_1 \psi_1(a) = b_1 \alpha(a), \qquad 
b_2 \psi_2(a) = b_2 \beta(a)
\end{equation*}
for $b_i \in \B_i, a \in \A.$
We set 
$\B_1 = \B_2 = \A$.
We put
${\cal H}_{\alpha,\beta} =\A
$ 
and
equip it with  Hilbert $C^*$-quad module structure
over $(\A;\A,\A)$
in the following way.
For $\xi = x, \xi' = x' \in 
{\cal H}_{\alpha, \beta}
 = \A,  a \in \A, 
z \in {\cal B}_1 =\A, w \in {\cal B}_2 =\A$,
define the right $\A$-module structure 
   and the right $\A$-valued inner product 
   $\langle \cdot \mid \cdot \rangle_\A$
   by
\begin{equation*}
\xi \cdot a = x a, 
\qquad 
\langle \xi \mid \xi' \rangle_\A = x^* x'.
\end{equation*}
Define the right actions
$\varphi_i$  of ${\cal B}_i$ 
with right $\B_i$-valued inner products 
$
\langle \cdot  \mid \cdot \rangle_{\B_i}
$
and
the left actions 
$\phi_i$ of ${\cal B}_i$
by setting
\begin{gather*}
\xi \varphi_1(z) = x {\alpha}(z), 
\qquad 
\xi \varphi_2(w) = x {\beta}(w), \\
\langle \xi \mid \xi' \rangle_{\B_1}
= {\alpha}^{-1}(x^*x'),
\qquad 
\langle \xi \mid \xi' \rangle_{\B_2}
= {\beta}^{-1}(x^*x'), \\
\phi_{1}(z) \xi = {\beta}({\alpha}(z))x, 
\qquad 
\phi_{2}(w) \xi =  {\alpha}({\beta}(w)) x.
\end{gather*}
It is straightforward to see  that 
${\cal H}_{\alpha, \beta}$
is a Hilbert $C^*$-quad module over $(\A;\A,\A)$  of strongly finite type.

\medskip

{\bf 2.}
We fix natural numbers $1 < N, M \in {\Bbb N}$.
Consider finite dimensional commutative $C^*$-algebras
$\A = {\Bbb C}, \ \B_1 = {\Bbb C}^N, \ \B_2 = {\Bbb C}^M.$
The right actions $\psi_i$ of $\A$
on $\B_i$ are naturally defined as 
right multiplications of 
${\Bbb C}$.
The algebras $\B_1, \B_2$ 
have the ordinary product structure and the inner product structure which we denote by 
$\langle \cdot \mid \cdot \rangle_N$ 
and
$\langle \cdot \mid \cdot \rangle_M$ 
respectively.
 Let us denote by
$\H_{M,N}$ the tensor product ${\Bbb C}^M \otimes {\Bbb C}^N$. 
Define the right actions 
$\varphi_i$ of $\B_i$  
with
$\B_i$-valued  right inner products 
$\langle \cdot \mid \cdot \rangle_{\B_i}
$
and
the left actions
$\phi_i$ of $\B_i$
on 
$\H_{M,N} = {\Bbb C}^M \otimes {\Bbb C}^N$
for $i=1,2$
by setting 
\begin{align*}
(\xi \otimes \eta )\varphi_1(z) 
 = \xi \otimes (\eta\cdot z), & \qquad
(\xi \otimes \eta )\varphi_2(w)
 = (\xi\cdot w) \otimes \eta,\\  
\langle \xi \otimes \eta \mid  \xi' \otimes \eta' \rangle_{\B_1}
& =\langle \xi \mid  \xi'  \rangle_{M} \eta^* \cdot \eta' \in \B_1, \\
\langle \xi \otimes \eta \mid  \xi' \otimes \eta' \rangle_{\B_2}
& =\xi^* \cdot \xi'  \langle \eta \mid  \eta'  \rangle_{N}\in \B_2, \\
\phi_1(z) (\xi \otimes \eta )
 = \xi \otimes (z \cdot \eta),& \qquad
\phi_2(w) (\xi \otimes \eta )
 = (w \cdot \xi) \otimes \eta
\end{align*}
for $z \in \B_1, \xi, \xi'\in {\Bbb C}^N, \, 
w \in \B_2, \eta, \eta'\in {\Bbb C}^M.$ 
Let $e_i, i=1,\dots,M$ 
and $f_k, k=1,\dots,N$
be the standard basis of 
${\Bbb C}^M $ and 
that of ${\Bbb C}^N$
respectively.
Put the finite bases
\begin{align*}
u_i & = e_i \otimes 1 \in \H_{M,N}, \qquad i=1,\dots,M,\\
v_k & = 1 \otimes f_k \in \H_{M,N}, \qquad i=1,\dots,N.
\end{align*}
It is straightforward to see  that 
$\H_{M,N}$ is a Hilbert $C^*$-quad module 
over $({\Bbb C}; {\Bbb C}^N,{\Bbb C}^M)$
of strongly finite type.

\medskip

{\bf 3.}
Let
$(\A,\rho,\eta, \Sigma^\rho,\Sigma^\eta,\kappa)$
be a $C^*$-textile dynamical system
which means that for $j \in \Sigma^\eta, l \in \Sigma^\rho$
endomorphisms $\eta_j, \rho_l$ of $\A$ 
are given with commutation relations
$\eta_j \circ \rho = \rho_k\circ \eta_i$ if $\kappa(l,j) = (i,k)$.
In \cite{MaPre20112},
a Hilbert $C^*$-quad module $\H^{\rho,\eta}_\kappa$ over $(\A;\B_1,\B_2)$
from 
$(\A,\rho,\eta, \Sigma^\rho,\Sigma^\eta,\kappa)$
is constructed (see \cite{MaPre20112} for its detail construction).
The two triplets 
$(\A,\rho, \Sigma^\rho)$
and 
$(\A,\eta, \Sigma^\eta)$
are $C^*$-symbolic dynamical systems
(\cite{MaCrelle}),
that yield $C^*$-algebras 
${\cal O}_\rho$
and
${\cal O}_\eta$
respectively.
The $C^*$-algebras 
$\B_1$ and $ \B_2$
are defined as the $C^*$-subalgebra of ${\cal O}_\eta$
generated by elements
$T_j y T_j^*, j \in 
\Sigma^\eta, y \in \A$
and 
that of ${\cal O}_\rho$
generated by
$S_k y S_k^*, k \in \Sigma^\rho, y \in \A$
respectively.
Define 
the maps
$\psi_i:\A \longrightarrow \B_i, i=1,2$
by
\begin{equation*}
\psi_1(a) = \sum_{j\in \Sigma^\eta} T_j a T_j^*,   \qquad
\psi_2(a) = \sum_{l\in \Sigma^\rho} S_l a S_l^* , \qquad a \in \A  
\end{equation*}
which yield the right actions of $\A$ on $\B_i, i=1,2$.
Define
the maps 
$\lambda_i: \B_i \longrightarrow \A, i=1,2$
 by 
\begin{equation*} 
\lambda_1(z) = \sum_{j\in \Sigma^\eta} T_j^* z T_j,  \qquad
\lambda_2(w) = \sum_{l\in \Sigma^\rho} S_l^* w S_l,
\qquad z \in \B_1, w \in \B_2.
\end{equation*}
Put
$e_j = T_j T_j^* \in \B_1, j\in \Sigma^\eta$.
Let
$z= \sum_{j \in \Sigma^\eta} T_j  z_j T_j^*$
be an element of
$\B_1$ for
$z_j \in \A$ with 
$T_j^* T_j z_j T_j^*T_j = z_j$.
As
$\lambda_1(e_j^* z) = \lambda_1(T_j z_j T_j^*) = z_j$,
one sees 
\begin{equation*}
z = \sum_{j \in \Sigma^\eta} T_j T_j^* T_j z_j T_j^*
= \sum_{j \in \Sigma^\eta} T_j T_j^* \psi_1( z_j)
= \sum_{j \in \Sigma^\eta} e_j \psi_1(\lambda_1(e_j^* z)).
\end{equation*}
We similarly have
by putting
$f_l = S_l S_l^*\in \B_2$,
\begin{equation*}
w = \sum_{l \in \Sigma^\rho}  f_l \psi_2(\lambda_2(f_l^* w))
\qquad
\text{ for }
w \in \B_2.
\end{equation*}
We see that
$\H^{\rho,\eta}_\kappa$
is a Hilbert $C^*$-quad module of strongly finite type.
In particular, two nonnegative commuting matrices $A, B$
with a specification $\kappa$ coming from the equality $AB = BA$ 
yield a $C^*$-textile dynamical system 
and hence a  Hilbert $C^*$-quad module of strongly finite type,
which are studied in  
\cite{MaPre2012}.
 
\section{Fock  Hilbert $C^*$-quad modules and creation operators}
In this section, 
we will 
construct a $C^*$-algebra from a Hilbert $C^*$-quad module $\H$ 
of general type
by using two kinds of creation operators on Fock space of Hilbert $C^*$-quad module.
We first
consider
relative tensor products of Hilbert $C^*$-quad modules
and then
 introduce Fock space of Hilbert $C^*$-quad modules
which is a two-dimensional analogue of Fock space of Hilbert $C^*$-bimodules.
We fix a Hilbert $C^*$-quad module
$\H$ over $(\A;\B_1,\B_2)$ of general type 
as in the preceding section.
The Hilbert $C^*$-quad module
$\H$ is originally a Hilbert $C^*$-right module over $\A$
with $\A$-valued inner product $\langle \cdot \mid \cdot \rangle_\A$.
It has two other structure of Hilbert $C^*$-bimodules,
the Hilbert $C^*$-bimodule
$(\phi_1,\H, \varphi_1)$ over $\B_1$  
and
the Hilbert $C^*$-bimodule
$(\phi_2,\H, \varphi_2)$ over $\B_2$
where $\phi_i$ is a left action of $\B_i$ on $\H$ and
$\varphi_i$ is a right action of $\B_i$ on $\H$ with
$\B_i$-valued right inner product 
$\langle \cdot \mid \cdot \rangle_{\B_i}$
for each $i=1,2$.
This situation is written as in the figure:
\begin{equation*}
\begin{CD} 
         @. \B_2 @.      \\
@.   @V{\phi_2}VV @. \\
\B_1 @>{\phi_1}>> \H @<{\varphi_1}<< \B_1 \\ 
@.   @A{\varphi_2}AA @. \\
         @. \B_2 @.      
\end{CD}
\end{equation*}
 We will define two kinds of relative tensor products
$$
\H \otimes_{\B_1} \H, 
\qquad
\H \otimes_{\B_2} \H
$$
as  Hilbert $C^*$-quad modules 
over $(\A;\B_1,\B_2)$.
The latter one should be written  
vertically as
$$
\begin{matrix}
\H \\
\otimes_{\B_2}\\
\H
\end{matrix}
$$
rather than horizontally
$\H 
\otimes_{\B_2}
\H.
$ 
The first relative tensor product
$\H\otimes_{\B_1}\H$
is defined as
the relative tensor product as Hilbert $C^*$-modules over $\B_1$,
where 
the left $\H$ is a right $\B_1$-module through $\varphi_1$
and   
the right $\H$ is a left $\B_1$-module through $\phi_1$.
It has a right $\B_1$-valued inner product 
and
a right $\B_2$-valued inner product defined by
\begin{align*}
\langle \xi\otimes_{\B_1} \zeta \mid \xi'\otimes_{\B_1} \zeta'\rangle_{\B_1}
& := 
\langle \zeta \mid \phi_{\B_1}(\langle \xi \mid \xi'\rangle_{\B_1})\zeta'\rangle_{\B_1}, \\
\langle \xi\otimes_{\B_1} \zeta \mid \xi'\otimes_{\B_1} \zeta'\rangle_{\B_2}
& := 
\langle \zeta \mid \phi_{\B_1}(\langle \xi \mid \xi'\rangle_{\B_1})\zeta'\rangle_{\B_2}
\end{align*}
 respectively.
It has two right actions, 
$\id \otimes\varphi_1$ from $\B_1$
and  
 $\id \otimes\varphi_2$ from $\B_2$.
It also has two left actions,
$\phi_1\otimes\id$ from $\B_1$
and $\phi_2\otimes\id$ from $\B_2$.
By these operations 
$\H\otimes_{\B_1}\H$
is a Hilbert $C^*$-bimodule over $\B_1$ 
as well as 
a Hilbert $C^*$-bimodule over $\B_2$.
It also has a  right $\A$-valued inner product defined by
\begin{equation*}
\langle \xi\otimes_{\B_1}\zeta \mid \xi'\otimes_{\B_1}\zeta'\rangle_\A
:= \lambda_1(\langle \xi\otimes_{\B_1}\zeta \mid \xi'\otimes_{\B_1}\zeta'\rangle_{\B_1})
(= \lambda_2(\langle \xi\otimes_{\B_1}\zeta \mid \xi'\otimes_{\B_1}\zeta'\rangle_{\B_2})),
\end{equation*}
a right $\A$-action $\id\otimes a$ for $a \in \A$
 and a left $\A$-action
$\phi_\A \otimes\id$.
By these structure
$\H\otimes_{\B_1}\H$
is a Hilbert $C^*$-quad module over $(\A;\B_1,\B_2)$.
\begin{equation*}
\begin{CD} 
         @. \B_2 @.      \\
@.   @V{\phi_2\otimes\id}VV @. \\
\B_1 @>{\phi_1\otimes\id}>> \H\otimes_{\B_1}\H @<{\id\otimes\varphi_1}<< \B_1 \\ 
@.   @A{\id\otimes\varphi_2}AA @. \\
         @. \B_2 @.      
\end{CD}
\end{equation*} 
We denote the above operations
$\phi_1\otimes\id, \phi_2\otimes\id,
\id\otimes\varphi_1, 
\id\otimes\varphi_2
$
still
by
$\phi_1, \phi_2,
\varphi_1, \varphi_2
$
respectively.
Similarly we consider the other  relative 
 tensor product
$\H\otimes_{\B_2}\H$
 defined by
the relative tensor product as Hilbert $C^*$-modules over $\B_2$,
where
the left $\H$ is a right $\B_2$-module through $\varphi_2$
and   
the right $\H$ is a left $\B_2$-module through $\phi_2$.
By a symmetric discussion to the above,
$\H\otimes_{\B_2}\H$
is a Hilbert $C^*$-quad module over $(\A;\B_1,\B_2)$.
The following lemma is routine.
\begin{lem}
Let
${\cal H}_i = \H, i=1,2,3$.
The correspondences
\begin{align*}
(\xi_1 \otimes_{\B_1} \xi_2) \otimes_{\B_2} \xi_3 \in ({\cal H}_1 \otimes_{\B_1} {\cal H}_2)\otimes_{\B_2} {\cal H}_3
& \longrightarrow 
\xi_1 \otimes_{\B_1} ( \xi_2 \otimes_{\B_2} \xi_3)  \in {\cal H}_1 \otimes_{\B_1}( {\cal H}_2 \otimes_{\B_2} {\cal H}_3), \\
(\xi_1 \otimes_{\B_2} \xi_2) \otimes_{\B_1} \xi_3 \in
({\cal H}_1 \otimes_{\B_2}  {\cal H}_2)\otimes_{\B_1} {\cal H}_3
& \longrightarrow
\xi_1 \otimes_{\B_2} (\xi_2 \otimes_{\B_1} \xi_3)
\in
{\cal H}_1 \otimes_{\B_2}( {\cal H}_2 \otimes_{\B_1} {\cal H}_3)
\end{align*}
yield  isomorphisms of Hilbert $C^*$-quad modules respectively.
\end{lem}
We write the isomorphism class of the former Hilbert $C^*$-quad modules as 
${\cal H}_1 \otimes_{\B_1}  {\cal H}_2\otimes_{\B_2} {\cal H}_3
$
and
that of the latter ones as
$
{\cal H}_1 \otimes_{\B_2}  {\cal H}_2\otimes_{\B_1} {\cal H}_3
$
respectively.

\medskip

Note that the direct sum
$\B_1 \oplus \B_2$ 
has a structure of  a pre Hilbert $C^*$-right module over $\A$
by the following operations:
For $b_1 \oplus b_2, b'_1 \oplus b'_2 \in \B_1 \oplus \B_2$
and $a \in \A$,
set
\begin{align*}
(b_1 \oplus b_2)  \psi_\A(a) 
& := b_1\psi_1(a) \oplus b_2  \psi_2(a) \in \B_1 \oplus \B_2, \\ 
\langle b_1 \oplus b_2  \mid b'_1 \oplus b'_2 \rangle_{\A} 
& := \lambda_1(b_1^* b'_1) + \lambda_2(b_2^* b'_2) \in \A.
\end{align*}
By \eqref{eqn:lambdabpsi}
the equality 
\begin{equation*}
\langle b_1 \oplus b_2  \mid (b'_1 \oplus b'_2) \psi_\A(a) \rangle_{\A} 
=
\langle b_1 \oplus b_2  \mid b'_1 \oplus b'_2 \rangle_\A\cdot a  
\end{equation*}
holds so that 
$\B_1 \oplus \B_2$ is a 
pre Hilbert $C^*$-right module over $\A$.
We denote by $F_0(\H)$ the completion of $\B_1 \oplus \B_2$
by the norm induced by the inner product
$\langle \cdot \mid \cdot \rangle_\A$.
It has
right $\B_i$-actions 
$\varphi_i$
and
left $\B_i$-action $\phi_i$
by
\begin{gather*}
(b_1 \oplus b_2) \varphi_1(z)
 =   b_1 z \oplus 0,  \qquad
(b_1 \oplus b_2) \varphi_2(w)
 = 0 \oplus  b_2 w,\\
\phi_1(z) (b_1 \oplus b_2) 
 =  z b_1 \oplus 0, \qquad
\phi_2(w) (b_1 \oplus b_2) 
 =  0 \oplus w b_2
\end{gather*}
for $b_1 \oplus b_2 \in \B_1 \oplus \B_2$,
$z \in \B_1, w \in \B_2$.


We denote the relative tensor product
$\H\otimes_{\B_i}\H$ 
and elements $\xi \otimes_{\B_i}\eta$
by 
 $\H\otimes_{i}\H$ 
and
$\xi \otimes_{i} \eta$
respectively for $i=1,2$.
Let us  define the Fock Hilbert $C^*$-quad module
as a two-dimensional analogue of the Fock space of Hilbert $C^*$-bimodules.
Put
$\Gamma_0 = \{ \emptyset \}$
and
$\Gamma_n = \{ (i_1,\dots,i_n)) \mid i_j = 1,2  \}, n=1,2,\dots.$
We set
\begin{align*}
F_1(\H) &  = \H, \\
F_2(\H) & = (\H\otimes_1 \HK) \oplus (\H\otimes_2 \H), \\
F_3(\H) & = (\H\otimes_1 \H \otimes_1\H) 
         \oplus (\H\otimes_1 \H \otimes_2 \H) \\
         & \oplus (\H\otimes_2 \H \otimes_1 \H) 
         \oplus (\H\otimes_2 \H \otimes_2 \H),\\
\cdots      & \cdots \\
F_n(\H) & 
= \oplus_{(i_1,\cdots, i_{n-1}) \in \Gamma_{n-1}}
\H \otimes_{i_1} \H\otimes_{i_2}\cdots \otimes_{i_{n-1}} \H \\
\cdots      & \cdots. 
 \end{align*}
as Hilbert $C^*$-bimodules over $\A$.
We will define the Fock Hilbert $C^*$-module
$F(\H)$ by setting
\begin{equation*}
F(\H) := \overline{\oplus_{n=0}^\infty F_n(\H)}
\end{equation*}
which is the completion of the algebraic direct sum
$\oplus_{n=0}^\infty F_n(\H)$
of the Hilbert $C^*$-right module over $\A$
under the norm
$\| \xi \|_\A $
on
$\oplus_{n=0}^\infty F_n(\H)$
induced by the $\A$-valued right inner product 
on $\oplus_{n=0}^\infty F_n(\H)$.
Then 
$F(\H)$ is a Hilbert $C^*$-right module over $\A$.
It has a natural
 left $\B_i$-action defined by $\phi_i$ 
for $i=1,2$.

\medskip

For $\xi\in \HK$
we define two operators 
\begin{align*}
s_\xi:& F_n(\H) \longrightarrow F_{n+1}(\H), \qquad n=0,1,2,\dots,\\
t_\xi:& F_n(\H) \longrightarrow F_{n+1}(\H), \qquad n=0,1,2,\dots
\end{align*}
by setting 
for $n=0$,
\begin{align*}
s_\xi(b_1 \oplus b_2) & = \xi\varphi_1(b_1), \qquad 
b_1 \oplus b_2 \in \B_1 \oplus \B_2,\\ 
t_\xi(b_1 \oplus b_2) & = \xi\varphi_2(b_2),  \qquad 
b_1 \oplus b_2 \in \B_1 \oplus \B_2,
\end{align*}
and
for $n=1,2,\dots$,
\begin{align*}
s_\xi 
(\xi_1 \otimes_{\pi_1}
\cdots  \otimes_{\pi_{n-1}}\xi_n)
& =
\xi\otimes_1 \xi_1 \otimes_{\pi_1}
\cdots  \otimes_{\pi_{n-1}}\xi_n, \\
t_\xi 
(\xi_1 \otimes_{\pi_1}
\cdots  \otimes_{\pi_{n-1}}\xi_n)
& =
\xi\otimes_2 \xi_1 \otimes_{\pi_1}
\cdots  \otimes_{\pi_{n-1}}\xi_n 
\end{align*}
for 
$\xi_1 \otimes_{\pi_1} \cdots \otimes_{\pi_{n-1}}\xi_n \in F_n(\H)$
with
$(\pi_1,\dots,\pi_{n-1})\in \Gamma_{n-1}$.
\begin{lem}
For $\xi\in \H$
the two operators 
\begin{align*}
s_\xi:& F_n(\H) \longrightarrow F_{n+1}(\H), \qquad n=0,1,2,\dots,\\
t_\xi:& F_n(\H) \longrightarrow F_{n+1}(\H), \qquad n=0,1,2,\dots
\end{align*}
are both right $\A$-module maps.
\end{lem}
\begin{pf}
We will show the assertion for $s_\xi$.
For $n=0$, 
we have for $b_1 \oplus b_2 \in \B_1 \oplus \B_2$ and $a \in \A$,
\begin{equation*}
s_\xi((b_1 \oplus b_2)\psi_\A(a)) 
 = \xi\varphi_1(b_1 \psi_1(a)) 
 = (\xi\varphi_1(b_1))a 
 = (s_\xi(b_1 \oplus b_2) )a. 
\end{equation*}
For $n=1,2,\dots $,
we have
\begin{align*}
 & s_\xi ((\xi_1 \otimes_{i_1}\xi_2 \otimes_{i_2}
 \cdots  \otimes_{i_{n-1}}\xi_n) a) \\
= & s_\xi ((\xi_1 \otimes_{i_1}\xi_2 \otimes_{i_2}
 \cdots  \otimes_{i_{n-1}}(\xi_n a)) 
=  \xi\otimes_{1} \xi_1 \otimes_{i_1}\xi_2 \otimes_{i_2}
 \cdots  \otimes_{i_{n-1}}(\xi_n a) \\
= &( \xi\otimes_{1} \xi_1 \otimes_{i_1}\xi_2 \otimes_{i_2}
 \cdots  \otimes_{i_{n-1}} \xi_n) a 
= [ s_\xi (\xi_1 \otimes_{i_1}\xi_2 \otimes_{i_2}
 \cdots  \otimes_{i_{n-1}} \xi_n)] a. 
 \end{align*}
\end{pf}
It is clear that 
the two operators 
$s_\xi, t_\xi$ yield bounded right $\A$-module maps 
on $F(\H)$ having its adjoints with respect to the $\A$-valued right
inner product on $F(\H)$.
The operators are still denoted by 
$s_\xi, t_\xi$ respectively. 
The adjoints of 
$s_\xi, t_\xi : F(\H) \longrightarrow F(\H)$
 with respect to the $\A$-valued  right inner product on $F(\H)$
 map $F_{n+1}(\H)$ to $F_n(\H)$, 
 $n=0,1,2,\dots$.
\begin{lem}
\begin{enumerate}
\renewcommand{\labelenumi}{(\roman{enumi})}
\item For $\xi, \xi'\in \H =F_1(\H)$, 
we have
\begin{equation*}
s_\xi^* \xi'  = \langle \xi \mid \xi' \rangle_{\B_1} \oplus 0 \in 
\B_1 \oplus \B_2, \qquad
t_\xi^* \xi'  = 0 \oplus \langle \xi \mid \xi' \rangle_{\B_2} \in 
\B_1 \oplus \B_2.
\end{equation*}
\item
For $\xi\in \H$
and
$\xi_1 \otimes_{i_1}\xi_2 \otimes_{i_2}
\cdots  \otimes_{i_n}\xi_{n+1} \in F_{n+1}(\H), n=1,2,\dots $, 
we have
\begin{align*}
& s_\xi^* 
(\xi_1 \otimes_{i_1}\xi_2 \otimes_{i_2}
\cdots  \otimes_{i_n}\xi_{n+1}) \\
=
&
{\begin{cases}
\phi_1(\langle \xi\mid  \xi_1\rangle_{\B_1})\xi_2 \otimes_{i_2}
\cdots  \otimes_{i_n}\xi_{n+1} 
   & \text{ if } i_1 = 1, \\
0  & \text{ if } i_1 = 2,
\end{cases}
} \\
& t_\xi^* 
(\xi_1 \otimes_{i_1}\xi_2 \otimes_{i_2}
\cdots  \otimes_{i_n}\xi_{n+1}) \\
=
&
{\begin{cases}
0  & \text{ if } i_1 = 1,\\
\phi_2(\langle \xi\mid  \xi_1\rangle_{\B_2})\xi_2 \otimes_{i_2}
\cdots  \otimes_{i_n}\xi_{n+1} 
   & \text{ if } i_1 = 2.
\end{cases}
}
\end{align*}
\end{enumerate}
\end{lem} 
\begin{pf}
We will show the assertions (i) and (ii) for $s_\xi^*$.
(i)
For $b_1 \oplus b_2 \in \B_1 \oplus \B_2$, 
we have
\begin{equation*}
\langle b_1 \oplus b_2 \mid s_\xi^* \xi' \rangle_\A
 = \langle \xi\varphi_1(b_1) \mid  \xi' \rangle_\A
  = \lambda_1(b_1^*\langle \xi \mid  \xi' \rangle_{\B_1}) 
   = \langle  b_1 \oplus b_2 \mid
\langle \xi \mid  \xi' \rangle_{\B_1} \oplus 0 \rangle_\A 
\end{equation*}
so that
$s_\xi^* \xi' = \langle \xi \mid \xi' \rangle_{\B_1} \oplus 0$.

(ii) 
For
$
\zeta_1 \otimes_{j_1}
\cdots
\otimes_{j_{n-1}}\zeta_n
\in F_n(\H)$ with $n=1,2,\dots$
we have 
\begin{align*}
& \langle
\zeta_1 \otimes_{j_1}
\cdots
\otimes_{j_{n-1}}\zeta_n
\mid
s_\xi^* 
(\xi_1 \otimes_{i_1}\xi_2 \otimes_{i_2}
\cdots  \otimes_{i_n}\xi_{n+1}) \rangle_\A \\
= & \langle
\xi \otimes_1 \zeta_1 \otimes_{j_1}
\cdots
\otimes_{j_{n-1}}\zeta_n
\mid
\xi_1 \otimes_{i_1}\xi_2 \otimes_{i_2}
\cdots  \otimes_{i_n}\xi_{n+1} \rangle_\A \\
= & \lambda_1(
\langle
\xi \otimes_1 \zeta_1 \otimes_{j_1}
\cdots
\otimes_{j_{n-1}}\zeta_n
\mid
\xi_1 \otimes_{i_1}\xi_2 \otimes_{i_2}
\cdots  \otimes_{i_n}\xi_{n+1} \rangle_{\B_1}) \\
= & 
{\begin{cases}
\langle
\zeta_1 \otimes_{j_1}
\cdots
\otimes_{j_{n-1}}\zeta_n
\mid
\phi_1(\langle \xi \mid \xi_1\rangle_{\B_1})\xi_2 \otimes_{i_2}
\cdots  \otimes_{i_n}\xi_{n+1} \rangle_\A 
  & \text{ if } i_1 = 1, \\
0 & \text{ if } i_1 = 2.
\end{cases}} 
\end{align*}
\end{pf}
Denote by
$\bar{\phi}_i $
the left actions of $\B_i, i=1,2$
on $F_n(\H)$ and hence on $F(\H)$ respectively.
They satisfy  the following equalities
\begin{align*}
\bar{\phi}_1(z)(b_1 \oplus b_2) = z b_1\oplus 0, 
& \qquad  
\bar{\phi}_2(w)(b_1 \oplus b_2) = 0 \oplus w b_2, \\
\bar{\phi}_1(z)(\xi_1 \otimes_{i_1}\xi_2 \otimes_{i_2}
\cdots  \otimes_{i_{n-1}}\xi_n)
& = 
(\phi_1(z)\xi_1) \otimes_{i_1}\xi_2 \otimes_{i_2}
\cdots  \otimes_{i_{n-1}}\xi_n,\\
\bar{\phi}_2(w)(\xi_1 \otimes_{i_1}\xi_2 \otimes_{i_2}
\cdots  \otimes_{i_{n-1}}\xi_n)
& = 
(\phi_2(w)\xi_1) \otimes_{i_1}\xi_2 \otimes_{i_2}
\cdots  \otimes_{i_{n-1}}\xi_n
\end{align*}
for $z \in \B_1, w \in B_2$, $b_1 \oplus b_2 \in \B_1 \oplus \B_2$
and $\xi_1 \otimes_{i_1}\xi_2 \otimes_{i_2}
\cdots  \otimes_{i_{n-1}}\xi_n \in F_n(\H)$.
More generally let us denote by
${\cal L}_\A(\H)$ and 
${\cal L}_\A(F(\H))$ 
the $C^*$-algebras of all bounded adjointable right $\A$-module maps on $\H$
and on $F(\H)$ with respect to 
their right $\A$-valued inner products respectively.
For $L \in {\cal L}_\A(\H)$,
define 
$\overline{L} \in {\cal L}_\A(F(\H))$
by 
\begin{align*}
\overline{L}(b_1\oplus b_2)  = 0  \quad \text{ for } 
& b_1\oplus b_2 \in \B_1 \oplus \B_2 \subset F_0(\H), \\
\overline{L}(\xi_1 \otimes_{i_1}\xi_2 \otimes_{i_2}
\cdots  \otimes_{i_{n-1}}\xi_n)
& = (L\xi_1) \otimes_{i_1}\xi_2 \otimes_{i_2}
\cdots  \otimes_{i_{n-1}}\xi_n
\end{align*} 
for
$\xi_1 \otimes_{i_1}\xi_2 \otimes_{i_2}
\cdots  \otimes_{i_{n-1}}\xi_n
\in F_n(\H).
$

\begin{lem}
Both the maps
$\bar{\phi}_i : \B_i \longrightarrow {\cal L}_\A(F(\H))$
for $i=1,2$
are faithful $*$-homomorphisms.
\end{lem} 
\begin{pf}
By assumption,
the $*$-homomorphisms
$\phi_i: \B_i \longrightarrow {\cal L}_\A(\H), i=1,2$
are faithful, 
so that the $*$-homomorphisms
$\bar{\phi}_i: \B_i \longrightarrow {\cal L}_\A(F(\H)), i=1,2$
are both faithful.
\end{pf}

\begin{lem}\label{lem:relations}
For $\xi, \zeta \in \H$, $z \in \B_1, w \in \B_2$,
$L \in {\cal L}_\A(\H)$ and $c,d \in {\Bbb C}$,
the following equalities hold on $F(\H)$:
\begin{gather}
s_{c \xi + d \zeta}  = c s_\xi + d s_\zeta, \qquad
t_{c \xi + d \zeta}  = c t_\xi + d t_\zeta,
\label{eqn:scxi}\\
s_{L \xi \varphi_1(z)}  = \overline{L} s_\xi \bar{\phi}_1(z),\qquad
t_{L \xi \varphi_2(w)}  = \overline{L} t_\xi \bar{\phi}_2(w),
\label{eqn:sLxi} \\
s_\zeta^* \overline{L} s_\xi 
  = \bar{\phi}_1( \langle \zeta \mid L \xi \rangle_{\B_1}), \qquad
t_\zeta^* \overline{L} t_\xi 
  = \bar{\phi}_2( \langle \zeta \mid L \xi \rangle_{\B_2}).
  \label{eqn:szetastar}
\end{gather}
\end{lem}
\begin{pf}
The equalities 
\eqref{eqn:scxi} are obvious.
We will show the equalities
\eqref{eqn:sLxi} and
\eqref{eqn:szetastar}
 for $s_\xi$.
We have for $b_1\oplus b_2 \in \B_1 \oplus \B_2$
\begin{align*}
s_{L \xi \varphi_1(z)}(b_1\oplus b_2) 
& =[L \xi \varphi_1(z)] \varphi_1(b_1) 
  =\overline{L}[ \xi \varphi_1(zb_1)] 
  =\overline{L}[ s_\xi (zb_1 \oplus 0)] \\
& =\overline{L}[ s_\xi [\bar{\phi}_1(z)(b_1 \oplus b_2)]] 
  =[\overline{L} s_\xi \bar{\phi}_1(z)](b_1 \oplus b_2).
\end{align*}
For 
$\xi_1 \otimes_{i_1}\xi_2 \otimes_{i_2}
\cdots  \otimes_{i_{n-1}}\xi_n \in F_n(\H)$,
$n=1,2,\dots $,
we have
\begin{align*}
s_{L\xi\varphi_1(z)}
(\xi_1 \otimes_{i_1}\xi_2 \otimes_{i_2}
\cdots  \otimes_{i_{n-1}}\xi_n)
& =(L\xi\varphi_1(z))
\otimes_1\xi_1 \otimes_{i_1}\xi_2 \otimes_{i_2}
\cdots  \otimes_{i_{n-1}}\xi_n \\
& =\overline{L} [(\xi\varphi_1(z))
\otimes_1\xi_1 \otimes_{i_1}\xi_2 \otimes_{i_2}
\cdots  \otimes_{i_{n-1}}\xi_n] \\
& =\overline{L} [\xi\otimes_1
(\phi_1(z)\xi_1) \otimes_{i_1}\xi_2 \otimes_{i_2}
\cdots  \otimes_{i_{n-1}}\xi_n] \\ 
& =\overline{L} [s_\xi((\phi_1(z)\xi_1) \otimes_{i_1}\xi_2 \otimes_{i_2}
\cdots  \otimes_{i_{n-1}}\xi_n)] \\ 
& =\overline{L} s_\xi \bar{\phi}_1(z)
[\xi_1 \otimes_{i_1}\xi_2 \otimes_{i_2}
\cdots  \otimes_{i_{n-1}}\xi_n]
\end{align*}
so that
$
s_{L\xi\varphi_1(z)}
=
\overline{L} s_\xi \bar{\phi}_1(z)
$
on $F_n(\H), n=0,1,\dots$.
Hence the equalities \eqref{eqn:sLxi} hold.

For $b_1 \oplus b_2 \in \B_1 \oplus \B_2$, 
we have 
\begin{align*}
 s_\zeta^* \overline{L}s_\xi(b_1 \oplus b_2) 
& = s_\zeta^*(L\xi\varphi_1(b_1)) 
  = \langle \zeta \mid L\xi\varphi_1(b_1) \rangle_{\B_1} \oplus 0\\
& = \langle \zeta \mid L\xi \rangle_{\B_1} b_1\oplus 0 
  = \bar{\phi}_1(\langle \zeta \mid L\xi \rangle_{\B_1}) (b_1 \oplus b_2). 
\end{align*}
For
$\xi_1\otimes_{i_1}\xi_2\otimes_{i_2}
\cdots
\otimes_{i_{n-1}}\xi_n
\in F_n(\H)$,
we have
\begin{align*}
s_\zeta^* \overline{L} s_\xi
(\xi_1\otimes_{i_1}\xi_2\otimes_{i_2}
\cdots
\otimes_{i_{n-1}}\xi_n)
 & =
s_\zeta^*\overline{L}(\xi \otimes_1 \xi_1\otimes_{i_1}\xi_2\otimes_{i_2}
\cdots
\otimes_{i_{n-1}}\xi_n)\\
&=
(\phi_1(\langle \zeta \mid L \xi \rangle_{\B_1})
\xi_1 )\otimes_{i_1}
\xi_2\otimes_{i_2}
\cdots
\otimes_{i_{n-1}}\xi_n\\
&=
\bar{\phi}_1(\langle \zeta \mid L \xi \rangle_{\B_1})
(\xi_1\otimes_{i_1}
\xi_2\otimes_{i_2}
\cdots
\otimes_{i_{n-1}}\xi_n
\end{align*}
so that
$s_\zeta^* \overline{L} s_\xi =\bar{\phi}_1(\langle \zeta \mid L \xi \rangle_{\B_1})$
on
$F_n(\H)$ for $n=0,1,2,\dots$.
Hence the equalities \eqref{eqn:szetastar} hold.
\end{pf}


The $C^*$-subalgebra of ${\cal L}_\A(F(\H))$
generated by 
the operators
$s_\xi, t_\xi$ for $\xi \in \H$
is denoted by
${\cal T}_{F(\H)}$ 
and is called the  Toeplitz quad module algebra
for $\H$.
\begin{lem}
The $C^*$-algebra
$\TFH$ contains the operators
$\bar{\phi}_1(z), \bar{\phi}_2(w)$ for $z \in \B_1, w \in \B_2$.
\end{lem}
\begin{pf}
By \eqref{eqn:szetastar}
in
the preceding lemma,
one sees 
$$
s_\zeta^* s_\xi 
  = \bar{\phi}_1( \langle \zeta \mid \xi \rangle_{\B_1}), \qquad
t_\zeta^*  t_\xi 
  = \bar{\phi}_2( \langle \zeta \mid \xi \rangle_{\B_2}),
  \qquad 
\zeta, \xi \in \H.
$$
Since $\H$ is a full $C^*$-quad module,
the inner products 
$\langle \zeta \mid \xi \rangle_{\B_1},
\langle \zeta \mid \xi \rangle_{\B_2}
$
for
$\zeta, \xi \in \H
$
generate
the $C^*$-algebras
$\B_1, \B_2$ respectively.
Hence 
$\bar{\phi}_1(\B_1), \bar{\phi}_2(\B_2)$
are contained in $\TFH$.
\end{pf}

\begin{lem}\label{lem:gauge}
There exists an action $\gamma$ of 
${\Bbb R}/{\Bbb Z} ={\Bbb T}$ on ${\cal T}_{F(\H)}$
such that
\begin{align*}
\gamma_r (s_\xi) & = e^{2 \pi \sqrt{-1} r } s_\xi,  \qquad 
\gamma_r (t_\xi)   = e^{2 \pi \sqrt{-1} r } t_\xi,  \qquad \xi \in \H,\\
\gamma_r (\bar{\phi}_1(z)) &=\bar{\phi}_1(z), \quad z \in \B_1,\qquad
\gamma_r (\bar{\phi}_1(w)) =\bar{\phi}_2(w), \quad w \in \B_2
\end{align*}
for $r \in {\Bbb R}/{\Bbb Z} = {\Bbb T}$.
\end{lem}
\begin{pf}
We will first define a one-parameter unitary group
$u_r, r \in {\Bbb R}/{\Bbb Z} ={\Bbb T}$
on
$F(\H)$ 
with respect to the 
right $\A$-valued inner product as in the following way.

For $n=0$ :
$ u_r : F_0(\H) \longrightarrow F_0(\H) 
$
is defined by
\begin{equation*}
u_r(b_1 \oplus b_2) = b_1 \oplus b_2 \quad \text{ for }
b_1 \oplus b_2 \in \B_1 \oplus \B_2.
\end{equation*}

For $n=1,2, \dots$ :
$ u_r : F_n(\H) \longrightarrow 
        F_n(\H) 
$
is defined by
\begin{equation*}
u_r(\xi_1\otimes_{i_1}\xi_2\otimes_{i_2}
\cdots
\otimes_{i_{n-1}}\xi_n)
 = e^{2 \pi \sqrt{-1}n r}\xi_1\otimes_{i_1}\xi_2\otimes_{i_2}
\cdots
\otimes_{i_{n-1}}\xi_n
\end{equation*}
 for
 $
\xi_1\otimes_{i_1}\xi_2\otimes_{i_2}
\cdots
\otimes_{i_{n-1}}\xi_n
 \in F_n(\H).
$
We therefore have a one-parameter unitary group
$u_r$ on $F(\H)$.
We then define an automorphism $\gamma_r$ on 
${\cal L}_\A(F(\H))$ 
for 
$r \in {\Bbb R}/ {\Bbb Z}$
by
\begin{equation*}
\gamma_r(T) = u_r T u_r^*
\quad
\text{ for }
T \in {\cal L}_\A(F(\H)), \, 
r \in {\Bbb R}/{\Bbb Z}.
\end{equation*}
It then follows that for $b_1 \oplus b_2 \in \B_1 \oplus \B_2$
\begin{equation*}
\gamma_r(s_\xi)(b_1\oplus b_2)
= u_r s_\xi u_r^* (b_1\oplus b_2) 
 = u_r (\xi(\varphi_1(b_1)) )
 = e^{2 \pi \sqrt{-1} r} s_\xi(b_1 \oplus b_2), 
\end{equation*}
and
for
$
\xi_1\otimes_{i_1}
\cdots
\otimes_{i_{n-1}}\xi_n
 \in F_n(\H),
n=1,2,\dots $
\begin{align*}
\gamma_r(s_\xi)
(\xi_1\otimes_{i_1}
\cdots
\otimes_{i_{n-1}}\xi_n)
& = e^{-2 \pi \sqrt{-1}n r} u_r 
(\xi\otimes_1 \xi_1\otimes_{i_1}
\cdots
\otimes_{i_{n-1}}\xi_n) \\
& = e^{2 \pi \sqrt{-1}r} 
s_\xi( \xi_1\otimes_{i_1}
\cdots
\otimes_{i_{n-1}}\xi_n). 
\end{align*}
Therefore we conclude that
$
\gamma_r(s_\xi) =
e^{2 \pi \sqrt{-1}r} s_\xi
$ 
on $F(\H)$
and similarly 
$
\gamma_r(t_\xi) =
e^{2 \pi \sqrt{-1}r} t_\xi
$ 
on $F(\H)$.
It is direct to see that 
\begin{equation*}
\gamma_r (\bar{\phi}_1(z)) =\bar{\phi}_1(z), 
\qquad
\gamma_r (\bar{\phi}_1(w)) =\bar{\phi}_2(w),
 \qquad
 \text{ for }
z \in \B_1, w \in \B_2.
\end{equation*}
It is also obvious that 
$\gamma_r(\TFH) =\TFH$ for
$r \in {\Bbb R}/{\Bbb Z}$.
\end{pf}

Denote by $J(\H)$
the $C^*$-subalgebra of
${\cal L}_\A(F(\H))$
generated 
by the elements
\begin{equation*}
{\cal L}_\A(\bigoplus_{n=0}^{\text{finite}} F_n(\H)).
\end{equation*}
The algebra $J(\H)$ is a closed two-sided  ideal of ${\cal L}_\A(F(\H))$.

\noindent
{\bf Definition.}
The $C^*$-algebra $\OFH$ associated to the Hilbert $C^*$-quad module $\H$
of general type is defined by the quotient $C^*$-algebra
of ${\cal T}_{F(\H)}$ by the ideal
${\cal T}_{F(\H)} \cap J(\H)$.

We denote by $[x]$ the quotient image of an element $x \in \TFH$
under the ideal 
${\cal T}_{F(\H)} \cap J(\H)$.
We set the elemets of $\OFH$  
\begin{equation*}
\S_\xi  = [s_\xi], \qquad
\T_\xi  = [t_\xi], \qquad
\varPhi_1(z)  = [\bar{\phi}_1(z)], \qquad
\varPhi_2(w)   = [\bar{\phi}_2(w)]
\end{equation*}
for $\xi \in \H$ and $z \in \B_1, w \in \B_2$.
By the preceding lemmas,
we have
\begin{prop}
The $C^*$-algebra $\OFH$ is generated by the family of operators
$\S_\xi, \T_\xi$ for $\xi \in \H$.
It contains the operators
$\varPhi_1(z), \varPhi_2(w)$ for 
$z \in \B_1, w \in \B_2$.
They satisfy the following equalities
\begin{gather}
\S_{c \xi + d \zeta}  = c \S_\xi + d \S_\zeta, \qquad
\T_{c \xi + d \zeta}  = c \T_\xi + d \T_\zeta, \label{eqn:prop1}
\\
\S_{\phi_1(z') \xi \varphi_1(z)} 
 = \varPhi_1(z') \S_\xi \varPhi_1(z),\qquad
\T_{\phi_1(z') \xi \varphi_2(w)}  
= \varPhi_1(z') \T_\xi \varPhi_2(w), \label{eqn:prop2}
 \\
\S_{\phi_2(w') \xi \varphi_1(z)} 
 = \varPhi_2(w') \S_\xi \varPhi_1(z),\qquad
\T_{\phi_2(w') \xi \varphi_2(w)}  
= \varPhi_2(w') \T_\xi \varPhi_2(w), \label{eqn:prop3}
 \\
 \S_\zeta^* \S_\xi 
  = \varPhi_1( \langle \zeta \mid  \xi \rangle_{\B_1}), \qquad
\T_\zeta^*  \T_\xi 
  = \varPhi_2( \langle \zeta \mid  \xi \rangle_{\B_2}) \label{eqn:prop4}
 \end{gather}
 for $\xi, \zeta \in \H$, $c,d \in {\Bbb C}$
and $z, z' \in \B_1, w, w' \in \B_2$.
\end{prop}

\section{The $C^*$-algebras of Hilbert $C^*$-quad modules of finite type}
In what follows, we assume that a Hilbert $C^*$-quad module $\H$
is of finite type.
In this section,
we will study the $C^*$-algebra 
$\OFH$
for a Hilbert $C^*$-quad module 
$\H$
of finite type.
Let $\{ u_1,\dots, u_M\}$ 
be a finite basis 
of $\H$ 
as a Hilbert $C^*$-right module over $\B_1$
and
$\{ v_1,\dots, v_N\}$
a finite basis 
of $\H$ 
as a Hilbert $C^*$-right module over $\B_2$.
Keep the notations as in the previous section.
We set
\begin{equation}
s_i  = s_{u_i} \quad \text{ for } i=1,\dots,M
\quad \text{ and }
\quad 
t_k  = t_{v_k} \quad \text{ for } k=1,\dots,N. 
\end{equation} 
By \eqref{eqn:basis} and Lemma \ref{lem:relations}
 we have for $\xi \in \H$
\begin{equation}
s_\xi  = \sum_{i=1}^M 
s_i \bar{\phi}_1(\langle u_i \mid \xi \rangle_{\B_1}), \qquad
t_\xi  = \sum_{k=1}^N 
t_{k} \bar{\phi}_2(\langle v_k \mid \xi \rangle_{\B_2}). \label{eqn:expan}
\end{equation}
Let $P_n$ be  the projection on $F(\H)$ 
onto $F_n(\H)$ for $ n=0,1,\dots$
so that 
$\sum_{n=1}^\infty P_n = 1$ on $F(\H)$.
\begin{lem}
For $\xi, \zeta\in \H$, 
we have 
\begin{enumerate}
\renewcommand{\labelenumi}{(\roman{enumi})}
\item  
$s_\xi^* t_\zeta P_n  = 0$ 
for $n=1,2,\dots$ and hence
$s_\xi^* t_\zeta  =s_\xi^* t_\zeta P_0$.
\item  
$t_\zeta^* s_\xi  P_n  = 0$
for $n=1,2,\dots$ and hence
$t_\zeta^* s_\xi   = t_\zeta^* s_\xi  P_0$.
\end{enumerate}
\end{lem}
\begin{pf}
(i)
For $n=1,2,\dots$,
we have 
\begin{equation*}
s_\xi^* t_\zeta ( \xi_1\otimes_{i_1}\cdots
\otimes_{i_{n-1}}\xi_n)
= 
s_\xi^*(\zeta \otimes_2 \xi_1\otimes_{i_1}
\cdots
\otimes_{i_{n-1}}\xi_n)
=0.
\end{equation*}

(ii) is smilarly shwon to (i).
\end{pf}
Define
two  projections on $F(\H)$ 
by
\begin{align*}
P_s
& = \text{The projection onto } 
\bigoplus_{n=0}^\infty
\sum_{(i_1,\cdots, i_{n}) \in \Gamma_{n}} 
\H \otimes_{1}\H\otimes_{i_1}\H\otimes_{i_2}\cdots
\otimes_{i_{n}}\H, \\
P_t
& = \text{The projection onto } 
\bigoplus_{n=0}^\infty
\sum_{(i_1,\cdots, i_{n})\in \Gamma_{n}}
\H \otimes_{2}\H\otimes_{i_1}\HK\otimes_{i_2}\cdots
\otimes_{i_{n}}\H.
\end{align*}

\begin{lem} 
Keep the above notations.
\begin{equation}
\sum_{i=1}^M s_i s_i^* = P_1 + P_s
\quad 
\text{ and }
\quad
\sum_{k=1}^N t_k t_k^* = P_1 + P_t. \label{eqn:stP}
\end{equation}
Hence
\begin{equation}
\sum_{i=1}^M s_i s_i^* 
+
\sum_{k=1}^N t_k t_k^*  + P_0 = 1_{F(\H)}  + P_1. \label{eqn:s+t}
\end{equation}
\end{lem}
\begin{pf}
For
$
\xi_1\otimes_{i_1}\xi_2\otimes_{i_2}
\cdots
\otimes_{i_{n-1}}\xi_n
\in F_n(\H)
$
with $2 \le n \in {\Bbb N}$,
we have
\begin{align*}
 &s_i s_i^*
  (\xi_1\otimes_{i_1}\xi_2\otimes_{i_2}
  \cdots
  \otimes_{i_{n-1}}\xi_n) \\
=&
{\begin{cases}
u_i\otimes_1
\phi_1(\langle u_i \mid \xi_1\rangle_\eta)\xi_2\otimes_{i_2}
\cdots
\otimes_{i_{n-1}}\xi_n)
   & \text{ if } i_1 = 1,\\
0  & \text{ if } i_1 = 2.
\end{cases} }
\end{align*}
As
$
u_i \otimes_1
\phi_1(\langle u_i \mid \xi_1\rangle_{\B_1})\xi_2
 = u_i 
\varphi_1(\langle u_i \mid \xi_1\rangle_{\B_1})
\otimes_1
\xi_2
$
and
$\sum_{i=1}^M u_i 
\varphi_1(\langle u_i \mid \xi_1\rangle_{\B_1})
= \xi_1$,
we have
\begin{equation*}
\sum_{i=1}^M s_i s_i^*
(\xi_1\otimes_{i_1}\xi_2\otimes_{i_2}
\cdots
\otimes_{i_{n-1}}\xi_n)
 =
\begin{cases}
\xi_1\otimes_{1}\xi_2\otimes_{i_2}
\cdots
\otimes_{i_{n-1}}\xi_n
   & \text{ if } i_1 = 1,\\
0  & \text{ if } i_1 = 2
\end{cases} 
\end{equation*}
and hence
\begin{equation*}
\sum_{i=1}^M s_i s_i^* |_{\oplus_{n=2}^{\infty}F_n(\H)}
= P_s|_{\oplus_{n=2}^{\infty}F_n(\H)}.
\end{equation*}
For $\xi \in F_1(\H) = \H$,
we have
$
s_i s_i^*
\xi
 = s_i(\langle u_i \mid \xi \rangle_{\B_1} \oplus 0)
= u_i \varphi_1(\langle u_i \mid \xi \rangle_{\B_1}) 
$
so that
\begin{equation*}
\sum_{i=1}^M s_i s_i^*
\xi
 = \sum_{i=1}^M
u_i \varphi_1(\langle u_i \mid \xi \rangle_{\B_1}) = \xi
\end{equation*}
and hence
\begin{equation*}
\sum_{i=1}^M s_i s_i^* |_{F_1(\H)}
= 1_{F_1(\H)}.
\end{equation*}
As
$ s_i s_i^*(b_1 \oplus b_2) = 0$ for 
$b_1 \oplus b_2 \in \B_1 \oplus \B_2$,
we have
\begin{equation*}
\sum_{i =1}^M s_i s_i^* |_{F_0(\H)}
= 0.
\end{equation*}
Therefore
we conclude that
\begin{equation*}
\sum_{i=1}^M s_i s_i^*
= P_s + P_1
\quad
\text{ and similarly }
\quad
\sum_{k=1}^N t_k t_k^* = P_t + P_1.
\end{equation*}
As 
$P_s + P_t + P_0 + P_1 = 1_{F(\H)},
$
one obtaines \eqref{eqn:s+t}.
\end{pf} 
We set the operators
\begin{align*}
\S_i  & = \S_{u_i}(= [s_i])  \quad \text{ for } i=1,\dots,M, 
\quad \text{ and } \\
\quad 
\T_k  & = \T_{v_k}(= [t_i])   \quad \text{ for } k=1,\dots,N 
\end{align*} 
in the $C^*$-algebra $\OFH$.
As two operators 
$
   \sum_{i=1}^M s_i s_i^*
$
and
$
\sum_{k=1}^N t_k t_k^*
$
are projections by \eqref{eqn:stP},
so are
$
\sum_{i=1}^M 
\S_i \S_i^*
$
and
$
 \sum_{k=1}^N 
 \T_k \T_k^*.
$
Since
$P_1 - P_0 \in J(\H)$,
the identity
\eqref{eqn:s+t} implies 
\begin{equation}
\sum_{i=1}^M 
\S_i \S_i^* +
 \sum_{k=1}^N 
 \T_k \T_k^* =1. \label{eqn:SiTk1}
\end{equation}
Therefore we have
\begin{thm}\label{thm:relations} 
Let $\H$ be a Hibert $C^*$-quad module over $(\A;\B_1,\B_2)$
of finite type
with finite basis $\{u_1,\dots,u_M\}$ as a right $\B_1$-module
and
$\{ v_1,\dots,v_N\}$  as a right $\B_2$-module.
Then we have
\begin{enumerate}
\renewcommand{\labelenumi}{(\roman{enumi})}
\item
The $C^*$-algebra $\OFH$
is generated  by the operators
$\S_1,\dots, \S_M, \T_1, \dots, \T_N$
and the elements
$\varPhi_1(z), \varPhi_2(w)$ for $z \in \B_1, w\in \B_2$.
\item
They satisfy the following operator relations:
\begin{align}
\sum_{i=1}^M \S_i \S_i^* +
\sum_{k=1}^N \T_k \T_k^* & =1, 
\qquad
\S_j^* \T_l =0, \label{eqn:tt1}
\\ 
\S_i^* \S_j  = \varPhi_1(\langle u_i \mid u_j \rangle_{\B_1}), 
& \qquad
\T_k^* \T_l  = \varPhi_2(\langle v_k \mid v_l \rangle_{\B_2}), \label{eqn:tt2} \\
\varPhi_1(z) \S_j 
= \sum_{i=1}^M \S_i 
\varPhi_1( \langle u_i \mid \phi_1(z) u_j \rangle_{\B_1}), 
& \qquad
\varPhi_1(z) \T_l 
= \sum_{k=1}^N \T_k 
\varPhi_2(\langle v_k \mid \phi_1(z) v_l \rangle_{\B_2}), \label{eqn:tt3}\\
\varPhi_2(w) \S_j 
= \sum_{i=1}^M \S_i 
\varPhi_1( \langle u_i \mid \phi_2(w) u_j \rangle_{\B_1}),
& \qquad
\varPhi_2(w) \T_l
= \sum_{k=1}^N \T_k 
\varPhi_2( \langle v_k \mid \phi_2(w) v_l \rangle_{\B_2}) \label{eqn:tt4}
\end{align}
for
$
z \in \B_1, w \in \B_2, i,j=1,\dots, M,
k,l=1,\dots,N$.
\item
There exists an action $\gamma$ of 
${\Bbb R}/{\Bbb Z} ={\Bbb T}$ on $\OFH$
such that
\begin{align*}
\gamma_r (\S_i) & = e^{2 \pi \sqrt{-1} r } \S_i,  
 \qquad 
\gamma_r (\T_k)   = e^{2 \pi \sqrt{-1} r } \T_k, \\
\gamma_r (\varPhi_1(z)) & =\varPhi_1(z), \qquad
\gamma_r (\varPhi_2(w)) =\varPhi_2(w) 
\end{align*}
for $r \in {\Bbb R}/{\Bbb Z} = {\Bbb T}$,
$i=1,\dots,M, k=1,\dots,N$,
and $ z \in \B_1, w \in \B_2$.
\end{enumerate}
\end{thm}
\begin{pf}
(i) The assertion comes from
the equalities
\eqref{eqn:expan}.

(ii) The first equality of \eqref{eqn:tt1} is \eqref{eqn:SiTk1}.
As the projection $P_0$ belongs to $J(H)$,
Lemma 4.1 ensures us the second equality of \eqref{eqn:tt1}.
The equalities \eqref{eqn:tt2} come from  
\eqref{eqn:prop4}.
For $z \in \B_1$ and $j=1,\dots,M,$
we have 
$ \phi_1(z)u_j = \sum_{i=1}^M u_i\varphi_1( \langle u_i \mid \phi_1(z)u_j\rangle_{\B_1}) 
$
so that
\begin{equation*}
\bar{\phi}_1(z) s_j = s_{\phi_1(z)u_j}
=  \sum_{i=1}^M s_{u_i} 
\bar{\phi}_1( \langle u_i \mid \phi_1(z)u_j\rangle_{\B_1}) 
\end{equation*}
which goes to the first equality of
 \eqref{eqn:tt3}.
The other equalities of 
 \eqref{eqn:tt3}  and \eqref{eqn:tt4}  
are similarly shown.

(iii) The assertion is direct from Lemma \ref{lem:gauge}.
\end{pf}
The action $\gamma$ of ${\Bbb T}$ on $\OFH$
defined in the above theorem (iii) is called 
the gauge action.

\section{The universal 
$C^*$-algebras associated with Hilbert $C^*$-quad modules}
In this section,
we will prove that the $C^*$-algebra 
$\OFH$ associated with a Hilbert $C^*$-quad module of finite type  
is the universal $C^*$-algebra subject to the operator relations
stated in Theorem \ref{thm:relations} (ii).
Throughout this section,
we fix a Hilbert $C^*$-quad module $\H$ 
over $(\A;\B_1,\B_2)$ of finite type
with finite basis 
$\{u_1,\dots,u_M\}$ as a right Hilbert $\B_1$-module and
$\{v_1,\dots,v_N\}$ as a right Hilbert $\B_2$-module 
as in the previous section.

Let $\PH$
be the universal $*$-algebra 
generated by operators
$S_1,\dots, S_M, T_1, \dots, T_N$
and elements
$z\in \B_1, w\in \B_2 $
subject to the relations: 
\begin{align}
\sum_{i=1}^M S_i S_i^* +
\sum_{k=1}^N T_k T_k^* & =1,
\qquad S_j^* T_l =0, \label{eqn:H1}\\ 
S_i^* S_j  = \langle u_i \mid u_j \rangle_{\B_1}, 
& \qquad
T_k^* T_l  = \langle v_k \mid v_l \rangle_{\B_2}, \label{eqn:H2} \\
z S_j = \sum_{i=1}^M S_i  \langle u_i \mid \phi_1(z) u_j \rangle_{\B_1}, 
& \qquad
z T_l = \sum_{k=1}^N T_k  \langle v_k \mid \phi_1(z) v_l \rangle_{\B_2},\label{eqn:H3} \\
w S_j = \sum_{i=1}^M S_i  \langle u_i \mid \phi_2(w) u_j \rangle_{\B_1},
& \qquad
w T_l = \sum_{k=1}^N T_k  \langle v_k \mid \phi_2(w) v_l \rangle_{\B_2} 
\label{eqn:H4}
\end{align}
for
$
z \in \B_1, w \in \B_2, i,j=1,\dots, M, \,
k,l=1,\dots,N$.
The above four relations 
\eqref{eqn:H1},\eqref{eqn:H2},\eqref{eqn:H3},\eqref{eqn:H4}
are called the relations $(\H)$.
In what follows, we fix 
operators
$S_1,\dots, S_M, T_1, \dots, T_N$
satisfying the relations $(\H)$.
\begin{lem}
The sums 
$
\sum_{i=1}^M S_i S_i^*
$
and
$
\sum_{k=1}^N T_k T_k^* 
$
are both projections.
\end{lem}
\begin{pf}
Put
$
P=\sum_{i=1}^M S_i S_i^*
$
and
$
Q =\sum_{k=1}^N T_k T_k^*. 
$
By the relations \eqref{eqn:H1}, 
one sees that
$0 \le P,Q \le 1$
and
$P+ Q =1, PQ =0$.
It is easy to see that
both $P $ and $Q$ are projections.
\end{pf}
\begin{lem}
\hspace{6cm}
\begin{enumerate}
\renewcommand{\labelenumi}{(\roman{enumi})}
\item
For $i,j =1,\dots,M$ and $z \in \B_1,w \in \B_2$ we have
\begin{equation*}
S_i^* z S_j  = \langle u_i \mid \phi_1(z) u_j \rangle_{\B_1},
\qquad  
S_i^* w S_j  = \langle u_i \mid \phi_2(w) u_j \rangle_{\B_1}. 
\end{equation*}
\item
For $k,l =1,\dots,N$ and $z \in \B_1,w \in \B_2$ we have
\begin{equation*}
T_k^* z T_l  = \langle v_k \mid \phi_1(z) v_l \rangle_{\B_2}, 
\qquad  
T_k^* w T_l  = \langle v_k \mid \phi_2(w) v_l \rangle_{\B_2}. 
\end{equation*}
\end{enumerate}
\end{lem}
\begin{pf}
(i)
By \eqref{eqn:H3}, we have
\begin{align*}
S_i^* z S_j 
& = \sum_{h=1}^M S_i^* S_h  \langle u_h \mid \phi_1(z) u_j \rangle_{\B_1} 
  = \sum_{h=1}^M \langle u_i \mid u_h \rangle_{\B_1}  
                 \langle u_h \mid \phi_1(z) u_j \rangle_{\B_1} \\
& =  \langle u_i \mid 
\sum_{h=1}^M u_h \langle u_h \mid \phi_1(z) u_j \rangle_{\B_1}  
\rangle_{\B_1}  
  =  \langle u_i \mid \phi_1(z) u_j \rangle_{\B_1}.  
\end{align*}
The other equality
$
S_i^* w S_j  = \langle u_i \mid \phi_2(w) u_j \rangle_{\B_1}
$
is similarly shown to the above equalities. 

(ii) is similar to (i).
\end{pf}

By the equalities
\eqref{eqn:uxieta},
\eqref{eqn:vxieta}
we have 
\begin{lem}
Keep the above notations. \hspace{3cm}
\begin{enumerate}
\renewcommand{\labelenumi}{(\roman{enumi})}
\item
For $w \in \B_2, j=1,\dots M$, 
the element 
$S_j^* w S_j$ belongs to $\A$ and the formula holds:
\begin{equation*}
\sum_{j=1}^M S_j^* \langle \xi \mid \eta \rangle_{\B_2}S_j =
\langle \xi \mid \eta \rangle_\A
\qquad 
\text{ for } \xi, \eta \in \H.
\end{equation*}
\item
For $z \in \B_1, l=1,\dots N$, 
the element 
$T_l^* z T_l$ belongs to $\A$ and the formula holds:
\begin{equation*}
\sum_{l=1}^N T_l^* \langle \xi \mid \eta \rangle_{\B_1}T_l =
\langle \xi \mid \eta \rangle_\A
\qquad 
\text{ for } \xi, \eta \in \H.
\end{equation*}
\end{enumerate}
\end{lem}

\begin{lem}
The following equalities 
for $z \in \B_1$ and $w \in \B_2$ hold:
\begin{enumerate}
\renewcommand{\labelenumi}{(\roman{enumi})}
\item
\begin{align}
z & = \sum_{i,j=1}^M  S_i \langle u_i \mid \phi_1(z)u_j \rangle_{\B_1} S_j^*
+
\sum_{k,l=1}^N  T_k \langle v_k \mid \phi_1(z)v_l \rangle_{\B_2} T_l^*,
\label{eqn:z}\\
w & = \sum_{i,j=1}^M  S_i \langle u_i \mid \phi_2(w)u_j \rangle_{\B_1} S_j^*
+
\sum_{k,l=1}^N  T_k \langle v_k \mid \phi_2(w)v_l \rangle_{\B_2} T_l^*. 
\label{eqn:w}
\end{align}
\item
\begin{align}
z^* & = \sum_{i,j=1}^M  
S_i \langle u_i \mid \phi_1(z)^* u_j \rangle_{\B_1} S_j^*
+
\sum_{k,l=1}^N  
T_k \langle v_k \mid \phi_1(z)^* v_l \rangle_{\B_2} T_l^*,
\label{eqn:z*}\\
w^* & = \sum_{i,j=1}^M  
S_i \langle u_i \mid \phi_2(w)^* u_j \rangle_{\B_1} S_j^*
+
\sum_{k,l=1}^N  T_k \langle v_k \mid \phi_2(w)^*
v_l \rangle_{\B_2} T_l^*. 
\label{eqn:w*}
\end{align}
\item
\begin{align*}
zw & = \sum_{i,j=1}^M  S_i \langle u_i 
\mid \phi_1(z) \phi_2(w)u_j \rangle_{\B_1} S_j^*
+
\sum_{k,l=1}^N  T_k \langle v_k 
\mid \phi_1(z) \phi_2(w) v_l \rangle_{\B_2} T_l^*,\\
wz & = \sum_{i,j=1}^M  S_i \langle u_i 
\mid \phi_2(w)\phi_1(z)u_j \rangle_{\B_1} S_j^*
+
\sum_{k,l=1}^N  T_k \langle v_k 
\mid \phi_2(w)\phi_1(z)v_l \rangle_{\B_2} T_l^*.
\end{align*}
\end{enumerate}
\end{lem}
\begin{pf}
(i)
By 
\eqref{eqn:H3} and 
\eqref{eqn:H4}, 
we have
\begin{equation*}
z S_j S_j^* 
 = \sum_{i=1}^M S_i \langle u_i \mid \phi_1(z) u_j \rangle_{\B_1}S_j^*, 
 \qquad
z T_l T_l^* 
 = \sum_{k=1}^N T_k \langle v_k \mid \phi_1(z) v_l \rangle_{\B_2}T_l^*   
\end{equation*}
so that by \eqref{eqn:H1} 
\begin{equation*}
z
 = \sum_{i,j =1}^M S_i \langle u_i \mid \phi_1(z) u_j \rangle_{\B_1}S_j^*
+  \sum_{k,l=1}^N T_k \langle v_k \mid \phi_1(z) v_l \rangle_{\B_2}T_l^*.   
\end{equation*}
Similarly we have \eqref{eqn:w}.

(ii)
All the adjoints
of $\phi_1(z), \phi_2(w)$
for $z \in \B_1, w \in \B_2$
by the three inner products
$
\langle \cdot \mid \cdot \rangle_{\B_1},
\langle \cdot \mid \cdot \rangle_{\B_2}$
and 
$
\langle \cdot \mid \cdot \rangle_{\A}$
on $\H$ coincide with $\phi_1(z^*), \phi_2(w^*)$
respectively.  
Hence the assertions are clear.

(iii)
By (i) we have
\begin{align*}
z w
=& ( \sum_{i,j =1}^M S_i \langle u_i \mid \phi_1(z) u_j \rangle_{\B_1}S_j^*
+   \sum_{k,l=1}^N T_k \langle v_k \mid \phi_1(z) v_l \rangle_{\B_2}T_l^*)\\
& \cdot 
(\sum_{g,h=1}^M  S_g \langle u_g \mid \phi_2(w)u_h \rangle_{\B_1} S_h^*
+
\sum_{m,n=1}^N  T_m \langle v_m \mid \phi_2(w)v_n \rangle_{\B_2} T_n^*).
\end{align*}
As $S_j^* T_m = T_l^* S_g= 0$
for any $j,g =1,\dots,M,\,
l,m=1,\dots,N$,
 it follows that
\begin{align*}
z w
=&  \sum_{i,j,g,h =1}^M 
S_i \langle u_i \mid \phi_1(z) u_j \rangle_{\B_1}S_j^*
S_g \langle u_g \mid \phi_2(w)u_h \rangle_{\B_1} S_h^* \\
+   
& \sum_{k,l,m,n=1}^N 
T_k \langle v_k \mid \phi_1(z) v_l \rangle_{\B_2}T_l^*
T_m \langle v_m \mid \phi_2(w)v_n \rangle_{\B_2} T_n^* \\
=&  \sum_{i,j,g,h =1}^M 
S_i \langle u_i \mid \phi_1(z) u_j \rangle_{\B_1}
\langle u_j \mid u_g \rangle_{\B_1} 
\langle u_g \mid \phi_2(w)u_h \rangle_{\B_1} S_h^* \\
+   
& \sum_{k,l,m,n=1}^N 
T_k \langle v_k \mid \phi_1(z) v_l \rangle_{\B_2}
\langle v_l \mid v_m \rangle_{\B_2}
 \langle v_m \mid \phi_2(w)v_n \rangle_{\B_2} T_n^* \\
=&  \sum_{i,g,h =1}^M 
S_i \langle u_i \mid \phi_1(z) 
\sum_{j=1}^M u_j \langle u_j \mid u_g \rangle_{\B_1}
\rangle_{\B_1}
\langle u_g \mid \phi_2(w)u_h \rangle_{\B_1} S_h^* \\
+ 
&  \sum_{k,m,n=1}^N 
T_k \langle v_k \mid \phi_1(z) 
\sum_{l=1}^N v_l \langle v_l \mid v_m \rangle_{\B_2}
\rangle_{\B_2}
 \langle v_m \mid \phi_2(w)v_n \rangle_{\B_2} T_n^* \\
=&  \sum_{i,h =1}^M 
S_i \langle u_i \mid \phi_1(z) 
\sum_{g=1}^M  u_g \rangle_{\B_1}
\langle u_g \mid \phi_2(w)u_h \rangle_{\B_1} S_h^* \\
+ 
&  \sum_{k,n=1}^N 
T_k \langle v_k \mid \phi_1(z) 
\sum_{m=1}^N v_m 
\rangle_{\B_2}
 \langle v_m \mid \phi_2(w)v_n \rangle_{\B_2} T_n^* \\
=&  \sum_{i,h =1}^M 
S_i \langle u_i \mid \phi_1(z) 
\phi_2(w)u_h \rangle_{\B_1} S_h^* 
+ 
 \sum_{k,n=1}^N 
T_k \langle v_k \mid \phi_1(z) 
 \phi_2(w)v_n \rangle_{\B_2} T_n^*. 
 \end{align*}
\end{pf}
\begin{lem}
Let $p(z,w)$ be a polynomial of elements of $\B_1$ and $\B_2$.
Then we have 
\begin{enumerate}
\renewcommand{\labelenumi}{(\roman{enumi})}
\item
$p(z,w) S_j = \sum_{i=1}^M S_i z_i$ for some $z_i \in \B_1$.
\item
$p(z,w) T_l = \sum_{k=1}^N T_k w_k$ for some $w_k \in \B_2$.
\end{enumerate}
\end{lem}
\begin{pf}
For $z \in \B_1, w \in \B_2$ and $i,j =1,\dots,M$,
by putting 
$z_{i,j} = \langle u_i \mid \phi_1(z) u_j \rangle_{\B_1} \in \B_1$
and
$w_{i,j} = \langle u_i \mid \phi_2(w) u_j \rangle_{\B_1} \in \B_1$,
 the relations \eqref{eqn:H3},\eqref{eqn:H4}
 imply
\begin{equation*}
z S_j  = \sum_{i=1}^M S_i z_{i,j}, \qquad 
w S_j  = \sum_{i=1}^M S_i w_{i,j}
\end{equation*}
so that the assertion of (i) holds.
(ii) is similarly shown to (i).
\end{pf}
\begin{lem}
Let $p(z,w)$ be a polynomial of elements of $\B_1$ and $\B_2$.
Then we have 
\begin{enumerate}
\renewcommand{\labelenumi}{(\roman{enumi})}
\item
$S_i^*p(z,w) S_j$ belongs to $\B_1$ for all $i,j = 1,\dots,M$.
\item
$T_k^*p(z,w) T_l$ belongs to $\B_2$ for all $k,l = 1,\dots,N$.
\item
$S_i^*p(z,w) T_l=0$ for all $i = 1,\dots,M, l = 1,\dots,N$.
\item
$T_k^*p(z,w) S_j =0$ for all $k = 1,\dots,N, j = 1,\dots,M$.
\end{enumerate}
\end{lem}
\begin{pf}
(i)
By the previous lemma, we know 
\begin{align*}
p(z,w) S_j 
& = \sum_{h=1}^M S_h z_h \qquad \text{ for some } z_h \in \B_1, \\
\intertext{ so that}
S_i^* p(z,w) S_j 
& = \sum_{h=1}^M S_i^* S_h z_h 
  = \sum_{h=1}^M \langle u_i \mid u_h \rangle_{\B_1} z_h. 
\end{align*}
As
$\langle u_i \mid u_h \rangle_{\B_1} z_h$ belongs to $\B_1$,
we see the assertion.
(ii) is similarly shown to (i).

(iii)
As $T_k^* S_i =0$, we have
\begin{equation*}
T_k^* p(z,w) S_j 
= \sum_{i=1}^M T_k^* S_i z_i =0. 
\end{equation*}
(iv) is similarly shown to (i).
\end{pf}
We set
\begin{equation}
S_{1,i} : = S_i,\quad  i=1,\dots, M,
\qquad
S_{2,k} : = T_k, \quad k=1,\dots, N.
\end{equation}
Put
$$
\Sigma_1 =\{ (1,i) \mid i=1,\dots,M \}, \qquad
\Sigma_2 =\{ (2,k) \mid k=1,\dots,N \}.
$$
\begin{lem}
Every element of $\PH$ can be written as a linear combination of elements 
of the  form
\begin{equation*}
S_{g_1,i_1}S_{g_2,i_2}\cdots S_{g_m,i_m} 
b 
S_{h_n,j_n}^* \cdots S_{h_2,j_2}^* S_{h_1,j_1}^*
\end{equation*}
for some
$(g_1,i_1),  (g_2, i_2), \dots, (g_m, i_m),
 (h_1,j_1),  (h_2, j_2), \dots, (h_n,j_n)
 \in \Sigma_1 \cup \Sigma_2$
where $b$ is a polynomial of elements of $\B_1$ and $\B_2$.
\end{lem}
\begin{pf}
The assertion follows from the preceding lemmas.
\end{pf}
By construction,
every representation of $\B_1$ and $\B_2$
on a Hilbert space $H$ together with operators 
$S_i,i=1,\dots, M,  T_k, k=1,\dots,N $ 
satisfying the relations $({\cal H})$
extends to a representation of $\PH$ on $B(H)$.
We will endow $\PH$ with the norm obtained by
taking the supremum of the norms in $B(H)$
over all such representations.
Note that this supremum is finite for every element of $\PH$
because of the inequalities 
$\| S_i \|, \| T_k \| \le 1$,
which come from \eqref{eqn:H1}.
The completion of the algebra $\PH$ under the norm
becomes a $C^*$-algebra denoted by $\OH$,
which is called the universal $C^*$-algebra subject to the relations $({\cal H})$.

Denote by $C^*(\phi_1(\B_1),\phi_2(\B_2))$
the $C^*$-subalgebra of ${\cal L}_\A(\H)$
generated by $\phi_1(\B_1)$ and $\phi_2(\B_2)$.
\begin{lem}
An element $L$ of the $C^*$-algebra
$C^*(\phi_1(\B_1),\phi_2(\B_2))$
is both a right $\B_1$-module map and a right $\B_2$-module map.
This means that the equalities
\begin{equation*}
[L \xi]\varphi_i(b_i)  = L [\xi \varphi_i(b_i)] \qquad
\text{ for } \xi \in \H, \, b_i \in \B_i
\end{equation*}
hold.
\end{lem}
\begin{pf}
Since both the operators
$\phi_1(z)$ for $z \in \B_1$
and
$\phi_2(w)$ for $w \in \B_2$
are right $\B_i$-module maps for $i=1,2$,
any element of the $*$-algebra algebraically generated by 
$\phi_1(\B_1)$ and 
$\phi_2(\B_2)$ is both
a right $\B_1$-module map and 
a right $\B_2$-module map.
Hence it is easy to see that
any element $L$ of the $C^*$-algebra
$C^*(\phi_1(\B_1),\phi_2(\B_2))$
is both a right $\B_1$-module map and a right $\B_2$-module map.
\end{pf}
Denote by 
$\B_\circ$
the $C^*$-subalgebra of $\OH$
generated by
$\B_1$ and $\B_2$.
\begin{lem}\label{lem:pi}
The correspondence
\begin{equation}
z,w \in \B_\circ \longrightarrow 
\phi_1(z), \phi_2(w) \in C^*(\phi_1(\B_1),\phi_2(\B_2))
\subset {\cal L}_\A(\H)
\end{equation}
gives rise to an isomorphism
from
$\B_\circ$ onto $C^*(\phi_1(\B_1),\phi_2(\B_2))$
as $C^*$-algebras.
\end{lem}
\begin{pf}
We note that by hypothesis
both the maps
\begin{align*}
\phi_1 & : z \in \B_1 \longrightarrow \phi_1(z) \in {\cal L}_\A(\H),\\
\phi_2 & : w \in \B_2 \longrightarrow \phi_2(w) \in {\cal L}_\A(\H)
\end{align*}
are injective.
Denote by
$P(\phi_1(\B_1),\phi_2(\B_2))$
the $*$-algebra on $\H$ algebraically generated by
$\phi_1(z), \phi_2(w)$
for
$z \in \B_1, w \in \B_2$. 
Define an operator 
$\pi(L) \in\OH$
for $L \in P(\phi_1(\B_1),\phi_2(\B_2))$
by
\begin{equation}
\pi(L) =
  \sum_{i,j=1}^M S_i \langle u_i \mid L u_j\rangle_{\B_1} S_j^*
+ \sum_{k,l=1}^N T_k \langle v_k \mid L v_l\rangle_{\B_2} T_l^*. 
\label{eqn:piL}
\end{equation}
Let
$P_\circ$
be the $*$-subalgebra of $\PH$ algebraically generated by $\B_1$ and $\B_2$.
Since
$\pi(\phi_1(z)) =z$ for $z \in \B_1$
and
$\pi(\phi_2(w)) = w$ for $w \in \B_2$
and by Lemma 5.4,
the map
$$
\pi : P(\phi_1(\B_1),\phi_2(\B_2)) \longrightarrow 
P_\circ \subset \B_\circ
$$
yields a $*$-homomorphism.
As $ S_i^*T_k = 0$
for $i=1,\dots,M, \, k=1,\dots,N,$ 
we have
\begin{equation*}
\| \pi(L) \|
= 
\Max \{ \| \sum_{i,j=1}^M S_i \langle u_i \mid L u_j\rangle_{\B_1} S_j^* \|, 
         \| \sum_{k,l=1}^N T_k \langle v_k \mid L v_l\rangle_{\B_2} T_l^* \|
\}.
\end{equation*}
We then have
\begin{equation*}
\| \sum_{i,j=1}^M S_i \langle u_i \mid L u_j\rangle_{\B_1} S_j^* \|
\le \sum_{i,j=1}^M \| \langle u_i \mid L u_j\rangle_{\B_1} \|
\le (\sum_{i,j=1}^M \| u_i \|_{\B_1} \| u_j \|_{\B_1}) \| L \| 
\end{equation*}
and similarly 
\begin{equation*}
 \| \sum_{k,l=1}^N T_k \langle v_k \mid L v_l\rangle_{\B_2} T_l^* \|
\le (\sum_{k,l=1}^N \| v_k \|_{\B_2} \| v_l \|_{\B_2}) \| L \|.
\end{equation*}
By putting
$
C = 
\Max \{ 
\sum_{i,j=1}^M \| u_i \|_{\B_1} \| u_j \|_{\B_1}, 
\sum_{k,l=1}^N \| v_k \|_{\B_2} \| v_l \|_{\B_2} 
\}$,
one has
$$
\| \pi(L) \| \le C \| L \| \qquad \text{ for all } 
L \in P(\phi_1(\B_1),\phi_2(\B_2)).
$$
Hence $\pi$  extends to 
the $C^*$-algebra 
$C^*(\phi_1(\B_1),\phi_2(\B_2))$
such that
$\pi(C^*(\phi_1(\B_1),\phi_2(\B_2))) = \B_\circ$.
The equality
\eqref{eqn:piL}
holds
for $L \in C^*(\phi_1(\B_1),\phi_2(\B_2))$.

We will next show that
$\pi: C^*(\phi_1(\B_1),\phi_2(\B_2)) \longrightarrow \B_\circ$ 
is injective. 
By \eqref{eqn:piL}, we have 
for
$L \in C^*(\phi_1(\B_1),\phi_2(\B_2))$
and
$h, h' =1,\dots,M$,
\begin{align*}
S_h^* \pi(L) S_{h'} 
& =
\sum_{i,j=1}^M S_h^* S_i \langle u_i \mid L u_j\rangle_{\B_1} S_j^* S_{h'} \\
& =
\sum_{i,j=1}^M \langle u_h \mid u_i \rangle_{\B_1}
               \langle u_i \mid L u_j\rangle_{\B_1} 
               \langle u_j \mid u_{h'} \rangle_{\B_1} \\
& =
\sum_{j=1}^M \langle u_h \mid 
\sum_{i=1}^M u_i \langle u_i \mid L u_j\rangle_{\B_1}
\rangle_{\B_1} \langle u_j \mid u_{h'} \rangle_{\B_1} \\
& =
\sum_{j=1}^M \langle u_h \mid L u_j\rangle_{\B_1}  
\langle u_j \mid u_{h'} \rangle_{\B_1} 
 =
 \langle  u_h  \mid L u_{h'} \rangle_{\B_1}. 
\end{align*}
Suppose that
$\pi(L) =0$
so that 
$
 \langle  u_h \mid L u_{h'} \rangle_{\B_1}=0.
$ 
Since
$$
L u_{h'} = \sum_{h=1}^M u_h
\langle  u_h \mid L u_{h'} \rangle_{\B_1},
$$
we see that
$
L u_{h'} = 0 
$ 
so that
$L =0$.
We thus conclude that
$\pi$ is injective and hence isomorphic.
\end{pf}
Denote by
$\phi_\circ : \B_\circ \longrightarrow C^*(\phi_1(\B_1), \phi_2(\B_2))$
the inverse $\pi^{-1}$ of the $*$-isomorphism $\pi$ giving in the proof of the above lemma
which satisfies
$$
\phi_\circ(z) = \phi_1(z) \quad \text{ for } z\in \B_1, \qquad
\phi_\circ(w) = \phi_2(w) \quad \text{ for } w\in \B_2.
$$
We put
${\cal F}_\H^0 = \B_\circ$.
For $n \in {\Bbb N}$, 
we denote by
$
{\cal F}_\H^n 
$
the closed linear span 
of elements of the form
\begin{equation*}
S_{g_1,i_1}S_{g_2,i_2}\cdots S_{g_n,i_n} 
b 
S_{h_n,j_n}^* \cdots S_{h_2,j_2}^* S_{h_1,j_1}^*
\end{equation*}
for some
$
(g_1,i_1),  (g_2, i_2), \dots, (g_n, i_n),
 (h_1,j_1),  (h_2, j_2), \dots, (h_n,j_n)
 \in \Sigma_1 \cup \Sigma_2
$
and $b \in \B_\circ$.
Let us denote by
${\cal F}_\H$ the $C^*$-subalgebra of $\OH$
generated by
$\cup_{n=0}^\infty {\cal F}_\H^n.$ 
By the relations \eqref{eqn:z} and \eqref{eqn:w},
we see the following 
\begin{lem}
For $x \in \B_\circ$, the following identity holds:
\begin{equation}
x  
= \sum_{i,j=1}^M  S_i \langle u_i \mid \phi_\circ(x)u_j \rangle_{\B_1} S_j^*
+
\sum_{k,l=1}^N  T_k \langle v_k \mid \phi_\circ(x)v_l \rangle_{\B_2} T_l^*.
\end{equation}
\end{lem}
Hence by putting for $b \in \B_\circ$
\begin{align*}
b_{1,ij} & = \langle u_i \mid \phi_\circ(b) u_j \rangle_{\B_1},
\qquad i,j = 1,\dots,M,\\
b_{2,kl} & = \langle v_k \mid \phi_\circ(b) v_l \rangle_{\B_2},
\qquad k,l = 1,\dots,N, 
\end{align*}
we have 
\begin{lem}
For $b \in \B_\circ$, the identity
\begin{align*}
& 
S_{g_1,i_1}S_{g_2,i_2}\cdots S_{g_n,i_n} 
b 
S_{h_n,j_n}^* \cdots S_{h_2,j_2}^* S_{h_1,j_1}^* \\ 
= 
& \sum_{i,j=1}^M  
S_{g_1,i_1}S_{g_2,i_2}\cdots S_{g_n,i_n}
S_{1,i} b_{1,ij} S_{1,j}^* 
S_{h_n,j_n}^* \cdots S_{h_2,j_2}^* S_{h_1,j_1}^* \\
+& 
\sum_{k,l=1}^N  
S_{g_1,i_1}S_{g_2,i_2}\cdots S_{g_n,i_n}
S_{2,k}  b_{2,kl}S_{2,l}^* 
S_{h_n,j_n}^* \cdots S_{h_2,j_2}^* S_{h_1,j_1}^*
\end{align*}
holds and induces an embedding of
${\cal F}_\H^n \hookrightarrow {\cal F}_\H^{n+1}$
for $n\in \Zp$. 
\end{lem}
\begin{lem}
The $C^*$-algebra ${\cal F}_\H$ 
is the inductive limit 
$\lim_{n\to\infty}{\cal F}_\H^n$
of the sequence of the inclusions: 
\begin{equation}
{\cal F}_\H^0    \hookrightarrow 
{\cal F}_\H^{1}  \hookrightarrow 
{\cal F}_\H^{2}  \hookrightarrow 
\cdots           \hookrightarrow 
{\cal F}_\H^n    \hookrightarrow 
{\cal F}_\H^{n+1}\hookrightarrow 
\cdots           \hookrightarrow 
{\cal F}_\H.
\end{equation}
\end{lem}
Let
$e^{2 \pi\sqrt{-1} r} \in  {\Bbb T}$
be a complex number of modulus one for
$r \in {\Bbb R}/{\Bbb Z}$.
The elements 
\begin{equation*}
e^{2 \pi\sqrt{-1} r} S_i, \, i=1,\dots,M, \quad
e^{2 \pi\sqrt{-1} r} T_k, \, k=1,\dots,N, \quad
z \in \B_1, \quad
w \in \B_2
\end{equation*}
in $\OH$ instead of
\begin{equation*}
S_i, \, i=1,\dots,M, \quad
T_k, \, k=1,\dots,N, \quad
z \in \B_1, \quad
w \in \B_2
\end{equation*}
satisfy the relations $(\H)$.
This implies the existence of an action on $\PH$ by automorphisms of 
the one-dimensional torus ${\Bbb T}$ that acts on the generators
by 
\begin{align*}
h_r(S_i) =e^{2 \pi\sqrt{-1} r} S_i,&  \qquad
h_r(T_k) =e^{2 \pi\sqrt{-1} r} T_k,  \\
h_r(z)= z, &  \qquad 
h_r(w) = w 
\end{align*}
for 
$i=1,\dots,M, k=1,\dots,N$,
$z \in \B_1, w \in \B_2$
and
$r \in {\Bbb R}/{\Bbb Z} ={\Bbb T}$. 
As the $C^*$-algebra $\OH$ 
has the largest norm on $\PH$,
the action $(h_r)_{r \in {\Bbb T}}$ on $\PH$ 
extends to an action of ${\Bbb T}$ on $\OH$,
still denote by $h$.
The formula
\begin{equation*}
a \in \OH \longrightarrow \int_{r \in {\Bbb T}} h_r(a) dr
\in \OH
\end{equation*} 
where $dr$ is the normalized Lebesgue measure on ${\Bbb T}$
defines a faithful conditional expectation 
denoted by
${\cal E}_\H$
from $\OH$ onto the fixed point algebra 
$(\OH)^h$.
The following lemma is routine.
\begin{lem}
$(\OH)^h = {\cal F}_\H.$
\end{lem}

\medskip

The $C^*$-algebra $\OH$ 
satisfies the following universal property.
Let ${\cal D}$ be a unital $C^*$-algebra and
$\Phi_1: \B_1 \longrightarrow {\cal D}$,
$\Phi_2: \B_2 \longrightarrow {\cal D}$
be $*$-homomorphisms 
such that
$\Phi_1(a) = \Phi_2(a)$ 
for 
$a \in \A$.
Assume that 
there exist elements 
$\widehat{S}_1,\dots, \widehat{S}_M, 
\widehat{T}_1, \dots, \widehat{T}_N$
in ${\cal D}$
satisfying the relations: 
\begin{align*}
\sum_{i=1}^M \widehat{S}_i \widehat{S}_i^* +
\sum_{k=1}^N \widehat{T}_k \widehat{T}_k^* &=1, 
\qquad
\widehat{S}_j^* \widehat{T}_l =0,
\\ 
\widehat{S}_i^* \widehat{S}_j  
= \Phi_1(\langle u_i \mid u_j \rangle_{\B_1}), 
& \qquad
\widehat{T}_k^* \widehat{T}_l  
= \Phi_2(\langle v_k \mid v_l \rangle_{\B_2}), \\
\Phi_1(z) \widehat{S}_j 
= \sum_{i=1}^M \widehat{S}_i  
\Phi_1(\langle u_i \mid \phi_1(z) u_j \rangle_{\B_1}), 
& \qquad
\Phi_1(z) \widehat{T}_l 
= \sum_{k=1}^N 
\widehat{T}_k  \Phi_2(\langle v_k \mid \phi_1(z) v_l \rangle_{\B_2}), \\
\Phi_2(w) \widehat{S}_j 
= \sum_{i=1}^M \widehat{S}_i  
\Phi_1(\langle u_i \mid \phi_2(w) u_j \rangle_{\B_1}),
& \qquad
\Phi_2(w) \widehat{T}_l = \sum_{k=1}^N 
\widehat{T}_k  
\Phi_2(\langle v_k \mid \phi_2(w) v_l \rangle_{\B_2}),
\end{align*}
for
$
z \in \B_1, w \in \B_2, i,j=1,\dots, M,
k,l=1,\dots,N$,
then there exists a unique $*$-homomorphism
$\Phi: \OH \longrightarrow {\cal D}$
such that 
\begin{equation*}
\Phi(S_i) = \widehat{S}_i,\quad
\Phi(T_k) = \widehat{T}_k,\quad
\Phi(z) = \Phi_1(z),\quad
\Phi(w) =  \Phi_2(w)
\end{equation*}
for   
$i=1,\dots, M, k=1,\dots,N$ 
and 
$z \in \B_1, w \in \B_2$.
We further assume that both the homomorphisms
$\Phi_i: \B_i \longrightarrow {\cal D}, i=1,2$
are injective.
We denote by 
$\Phi_\circ :\B_\circ \longrightarrow {\cal D}$
the restriction of $\Phi$ to the subalgebra 
$\B_\circ$.
 Let us denote by
 $\widehat{\cal O}_\H$ 
 the $C^*$-subalgebra of ${\cal D}$
 generated by
$
 \widehat{S}_i,
\widehat{T}_k, i=1,\dots, M, k=1,\dots,N
$
and
$
\Phi_1(z), \Phi_2(w) 
$
for
$z \in \B_1$, $w \in \B_2$.
\begin{lem}\label{lem:phi0}
Keep the above situation.
The $*$-homomorphism
$\Phi_\circ: \B_\circ \longrightarrow {\cal D}$
is injective.
\end{lem}
\begin{pf}
Since the correspondence in Lemma 5.9 
\begin{equation*}
\phi_\circ: z,w,\in \B_\circ \longrightarrow 
\phi_1(z), \phi_2(w) \in C^*(\phi_1(\B_1),\phi_2(\B_2))
\end{equation*}
yields an isomorphism of $C^*$-algebras,
it suffices to prove that the correspondence 
\begin{equation*}
\phi_1(z), \phi_2(w) \in C^*(\phi_1(\B_1),\phi_2(\B_2))
\longrightarrow
\Phi_1(z), \Phi_2(w) \in {\cal D}
\end{equation*}
 yields an isomorphism.
Let
$\widehat{\B}_\circ $  
be the $C^*$-subalgebra of $\widehat{\cal O}_\H$
generated by elements
$
\Phi_1(z), \Phi_2(w) \in {\cal A}
$
for $z\in \B_1, w \in  \B_2$.
Define  an element
$\widehat{\pi}(L)$
of
${\cal D}$
for 
$
L \in C^*(\phi_1(\B_1),\phi_2(\B_2))
$
by setting
\begin{equation}
\widehat{\pi}(L) =
  \sum_{i,j=1}^M \widehat{S}_i 
  \Phi_1(\langle u_i \mid L u_j\rangle_{\B_1}) 
  \widehat{S}_j^*
+ \sum_{k,l=1}^N \widehat{T}_k 
\Phi_2(\langle v_k \mid L v_l\rangle_{\B_2}) 
\widehat{T}_l^* \in {\cal D}.
\end{equation}
As in the proof of Lemma 5.9,
one sees that 
$\widehat{\pi}$ gives rise to a $*$-homomorphism
from
$C^*(\phi_1(\B_1),\phi_2(\B_2))$ into ${\cal D}$.
Since
\begin{equation*}
\widehat{\pi}(\phi_1(z)) =
  \sum_{i,j=1}^M \widehat{S}_i 
  \Phi_1(\langle u_i \mid \phi_1(z) u_j\rangle_{\B_1}) 
  \widehat{S}_j^*
+ \sum_{k,l=1}^N \widehat{T}_k 
\Phi_2(\langle v_k \mid \phi_1(z) v_l\rangle_{\B_2}) 
\widehat{T}_l^* 
= \Phi_1(z),
\end{equation*}
and similarly
$\widehat{\pi}(\phi_2(w)) = \Phi_2(w)$,
it is enough to show that
$\widehat{\pi}$ 
is injective. 
Suppose that
$\widehat{\pi}(L) =0$ for some 
 $L \in C^*(\phi_1(\B_1),\phi_2(\B_2))$.
By following the proof of Lemma \ref{lem:pi},
one sees that
$
\widehat{S}_h^* \widehat{\pi}(L) \widehat{S}_{h'}
 =
\Phi_1( \langle  u_h  \mid L u_{h'} \rangle_{\B_1}) 
$
for all $h, h' =1,\dots, M$. 
Hence the condition
$\widehat{\pi}(L) =0$ 
implies 
$\Phi_1( \langle  u_h  \mid L u_{h'} \rangle_{\B_1})=0$.
Since
$\Phi_1$ is injective,
we have
$\langle  u_h  \mid L u_{h'} \rangle_{\B_1}=0$
for all $h, h' =1,\dots, M$. 
As $L$ is a right $\B_1$-module map,
we have
for $\xi \in \H$,
$$
L\xi = \sum_{h' =1}^M L ( u_{h'}\langle u_{h'} \mid \xi \rangle_{\B_1} )
     = \sum_{h' =1}^M (L u_{h'} )\langle u_{h'} \mid \xi \rangle_{\B_1} 
     =0
$$
so that $L =0$.
Therefore 
$\widehat{\pi}: C^*(\phi_1(\B_1),\phi_2(\B_2)) \longrightarrow {\cal D}$ 
is injective. 
Hence the composition
$$
\Phi_\circ: \widehat{\pi} \circ \phi_\circ:
\B_\circ \overset{\phi_\circ}{\longrightarrow} 
C^*(\phi_1(\B_1),\phi_2(\B_2)) 
\overset{\widehat{\pi}}{\longrightarrow} {\cal D}
$$ 
is injective. 
\end{pf}
We set
\begin{equation}
\widehat{S}_{1,i} : = \widehat{S}_i,\quad  i=1,\dots, M,
\qquad
\widehat{S}_{2,k} : = \widehat{T}_k, \quad k=1,\dots, N.
\end{equation}
We put
$
\widehat{{\cal F}}_\H^0
=\B_\circ. 
$
For $n \in {\Bbb N}$, 
let 
$
\widehat{{\cal F}}_\H^n 
$
be the closed linear span 
in the $C^*$-algebra $\widehat{\cal O}_\H$
of elements of the form
$$
\widehat{S}_{g_1,i_1} \widehat{S}_{g_2,i_2}\cdots \widehat{S}_{g_n,i_n} 
\Phi_\circ(b)
\widehat{S}_{h_n,j_n}^* \cdots \widehat{S}_{h_2,j_2}^* \widehat{S}_{h_1,j_1}^*
$$
for
$
(g_1,i_1),  (g_2, i_2), \dots, (g_n, i_n),
 (h_1,j_1),  (h_2, j_2), \dots, (h_n,j_n)
 \in \Sigma_1 \cup \Sigma_2
$
and $b \in \B_\circ$.
Similarly to the subalgebras
$
{\cal F}_\H^n, n \in \Zp,
$
of $\OH$,
one knows that 
the closed linear span
$
\widehat{{\cal F}}_\H^n 
$
is a $C^*$-algebra and 
naturally regarded as a subalgebra of
$
\widehat{{\cal F}}_\H^{n+1} 
$
for each $n \in \Zp$.
Let us denote by
$\widehat{{\cal F}}_\H$ the $C^*$-subalgebra of 
$\widehat{{\cal O}}_\H$
generated by
$\cup_{n=0}^\infty \widehat{{\cal F}}_\H^n.$ 
Then the $C^*$-algebra $\widehat{{\cal F}}_\H$ 
is the inductive limit 
$\lim_{n\to\infty}\widehat{{\cal F}}_\H^n$
of the sequence of the inclusions :
\begin{equation}
\widehat{{\cal F}}_\H^0    \hookrightarrow 
\widehat{{\cal F}}_\H^{1}  \hookrightarrow 
\widehat{{\cal F}}_\H^{2}  \hookrightarrow 
\cdots           \hookrightarrow 
\widehat{{\cal F}}_\H^n    \hookrightarrow 
\widehat{{\cal F}}_\H^{n+1}\hookrightarrow 
\cdots           \hookrightarrow 
\widehat{{\cal F}}_\H. \label{eqn:AF}
\end{equation}
\begin{lem}
Suppose that both the $*$-homomorphisms
$\Phi_i: \B_i \longrightarrow \widehat{{\cal O}}_\H, i=1,2$
are injective.
Then the restriction of
$\Phi$ to the subalgebra ${\cal F}_\H$ yields a $*$-isomorphism
$\Phi|_{{\cal F}_\H}:{\cal F}_\H \longrightarrow \widehat{{\cal F}}_\H$.
\end{lem}
\begin{pf}
By the universality of $\OH$,
the restriction of
$\Phi$ to ${\cal F}_\H$ yields a surjective $*$-homomorphism
$\Phi|_{{\cal F}_\H}:{\cal F}_\H \longrightarrow \widehat{{\cal F}}_\H$.
It suffices to show that
$\Phi|_{{\cal F}_\H}$ is injective.
Suppose that $\Ker(\Phi|_{{\cal F}_\H}) \ne \{0 \}$
and put
$I =\Ker(\Phi|_{{\cal F}_\H})$.
Since
$\Phi_{\cal F}({\cal F}_\H^n) =\widehat{{\cal F}}_\H^n$
and
${\cal F}_\H = \lim_{n\to \infty}{\cal F}_\H^n$,
there exists $n \in \Zp$
such that
$I \cap {\cal F}_\H^n \ne 0$.
Let us denote by
$\Sigma^n$ the set of $n$-tuples of
$\Sigma_1 \cup \Sigma_2$:
$$
\Sigma^n =\{ (\mu_1, \dots,\mu_n) \mid 
\mu_1, \dots,\mu_n \in \Sigma_1 \cup \Sigma_2 \}.
$$
For 
$\mu = (\mu_1,\dots,\mu_n) \in \Sigma^n$,
denote by
$S_\mu$ the operator
$$
S_\mu = S_{\mu_1}S_{\mu_2}\cdots S_{\mu_n}
$$ 
where
$$
S_{\mu_m} = 
\begin{cases}
S_{1,i} = S_i  \quad & \text{ if } \mu_m = (1,i) \in \Sigma_1,\\ 
S_{2,k} = T_k  \quad & \text{ if } \mu_m = (2,k) \in \Sigma_2, \, m=1,\dots,n.
\end{cases}
$$
Any element of 
${\cal F}_\H^n $
is of the form
\begin{equation*}
\sum_{\mu,\nu\in \Sigma^n} S_\mu b_{\mu,\nu}S_\nu^*
\quad
\text{ for some }
b_{\mu,\nu} \in \B_\circ.
\end{equation*}
Hence one may find a nonzero 
element
$\sum_{\mu,\nu\in \Sigma^n} S_\mu b_{\mu,\nu}S_\nu^* 
\in I \cap {\cal F}_\H^n$.
Since
$\sum_{i=1}^M S_i S_i^* + \sum_{k=1}^N T_k T_k^*=1$,
the equality 
$\sum_{\mu \in \Sigma^n} S_\mu S_\mu^*=1$
hols.
For any $\omega, \gamma \in \Sigma^n$, 
one then sees 
$$
0 \ne S_\omega^* ( \sum_{\mu,\nu\in \Sigma^n} S_\mu b_{\mu,\nu}S_\nu^*) S_\gamma 
\in I \cap {\cal F}_\H^n.
$$
As 
$S_i^* T_k =0$
and
$S_i^*S_j = \langle u_i \mid u_j \rangle_{\B_1}$, 
$T_k^*T_l = \langle v_k \mid v_l \rangle_{\B_2}$
for
$i,j=1,\dots,M, \, k,l=1,\dots,N$,
the element
$S_\omega^* ( \sum_{\mu,\nu\in \Sigma^n} S_\mu b_{\mu,\nu}S_\nu^*) S_\gamma
$
belongs to
$I \cap \B_\circ$.
By the preceding lemma,
the homomorphism
$\Phi_\circ: \B_\circ \longrightarrow \widehat{\cal O}_\H$
is injective,
so that we have
$\Phi_\circ(S_\omega^* ( \sum_{\mu,\nu\in \Sigma^n} S_\mu b_{\mu,\nu}S_\nu^*) S_\gamma) \ne 0$
a contradiction.
Therefore we conclude that 
$\Phi|_{{\cal F}_\H}:{\cal F}_\H \longrightarrow \widehat{{\cal F}}_\H$
is injective and hence isomorphic.
\end{pf}

The following theorem is one of the main results of the paper.
\begin{thm}\label{thm:main1}
Let ${\cal D}$ be a unital $C^*$-algebra.
Suppose that there exist $*$-homomorphisms
$\Phi_1: \B_1 \longrightarrow {\cal D}$,
$\Phi_2: \B_2 \longrightarrow {\cal D}$
such that
$\Phi_1(a) = \Phi_2(a)$
for $a \in \A$
and 
there exist elements 
$\widehat{S}_1,\dots, \widehat{S}_M, 
\widehat{T}_1, \dots, \widehat{T}_N$
in ${\cal D}$
satisfying the relations: 
\begin{align*}
\sum_{i=1}^M \widehat{S}_i \widehat{S}_i^* +
\sum_{k=1}^N \widehat{T}_k \widehat{T}_k^* & =1,
\qquad 
S_j^* T_l =0, \\ 
\widehat{S}_i^* \widehat{S}_j  
= \Phi_1(\langle u_i \mid u_j \rangle_{\B_1}), 
& \qquad
\widehat{T}_k^* \widehat{T}_l  
= \Phi_2(\langle v_k \mid v_l \rangle_{\B_2}), \\
\Phi_1(z) \widehat{S}_j 
= \sum_{i=1}^M \widehat{S}_i  
\Phi_1(\langle u_i \mid \phi_1(z) u_j \rangle_{\B_1}), 
& \qquad
\Phi_1(z) \widehat{T}_l 
= \sum_{k=1}^N 
\widehat{T}_k  \Phi_2(\langle v_k \mid \phi_1(z) v_l \rangle_{\B_2}), \\
\Phi_2(w) \widehat{S}_j 
= \sum_{i=1}^M \widehat{S}_i  
\Phi_1(\langle u_i \mid \phi_2(w) u_j \rangle_{\B_1}),
& \qquad
\Phi_2(w) \widehat{T}_l = \sum_{k=1}^N 
\widehat{T}_k  
\Phi_2(\langle v_k \mid \phi_2(w) v_l \rangle_{\B_2}),
\end{align*}
for
$
z \in \B_1, w \in \B_2, i,j=1,\dots, M,
k,l=1,\dots,N$.
 Let us denote by
 $\widehat{\cal O}_\H$ 
 the $C^*$-subalgebra of ${\cal D}$
 generated by
$
 \widehat{S}_i,
\widehat{T}_k, i=1,\dots, M, k=1,\dots,N
$
and
$
\Phi_1(z), \Phi_2(w), 
$
for
$z \in \B_1, w \in \B_2$.
We further assume that 
the algebra
 $\widehat{\cal O}_\H$ 
 admits a gauge action.
If both the $*$-homomorphisms
$\Phi_i: \B_i \longrightarrow {\cal A}, i=1,2$
are injective,
then
there exists a $*$-isomomorphism
$\Phi: \OH \longrightarrow \widehat{\cal O}_\H$
satisfying
\begin{equation}
\Phi(S_i) = \widehat{S}_i,\quad
\Phi(T_k) = \widehat{T}_k,\quad
\Phi(z) = \Phi_1(z),\quad
\Phi(w) =  \Phi_2(w)
\end{equation}
for   
$i=1,\dots, M, k=1,\dots,N$ 
and 
$z \in \B_1, w \in \B_2$.
\end{thm}
\begin{pf}
By assumption,
 $\widehat{\cal O}_\H$ 
 admits a gauge action,
which we dnote by $\widehat{h}$.
Let us denote 
by $(\widehat{\cal O}_\H)^{\widehat{h}}$
the fixed point algebra of 
$\widehat{\cal O}_\H$ 
under the gauge action $\widehat{h}$
and
by $\widehat{\cal F}_\H$
the $C^*$-subalgebra 
of
$\widehat{\cal O}_\H$ 
defined by the inductive limit 
\eqref{eqn:AF}.
Then it is routine to check that
$(\widehat{\cal O}_\H)^{\widehat{h}}$ 
is canonically $*$-isomorphic to 
$\widehat{\cal F}_\H$.
There exists a conditional expectation
$$
\widehat{{\cal E}}_\H: \widehat{\cal O}_\H \longrightarrow 
\widehat{\cal F}_\H  
$$
defined by 
$$
\widehat{{\cal E}}_\H(x) = \int_{r \in {\Bbb T}} \widehat{h}_r(x) dr
\quad
\text{ for }
x \in \widehat{\cal O}_\H.
$$
By the universality of the algebra 
$\OH$ there exists a surjective $*$-homomorphism
$\Phi$ from $\OH$ to
$\widehat{\cal O}_\H$ such that 
$$
\Phi(S_i) = \widehat{S}_i,\quad 
\Phi(T_k) = \widehat{T}_k, \qquad
\Phi(z) = \Phi_1(z), \quad
\Phi(w) = \Phi_2(w)
$$
for
$
i,j=1,\dots, M,
k,l=1,\dots,N,
\, z \in \B_1, w \in \B_2$.
Then 
$\Phi({\cal F}_\H)=
\widehat{\cal F}_\H$
and the following diagram:
$$
\begin{CD}
{\cal O}_\H @>\Phi>>\widehat{O}_\H \\
@V{\cal E}_\H VV @VV \widehat{\cal E}_\H V \\ 
{\cal F}_\H @>\Phi|_{{\cal F}_\H}>> \widehat{F}_\H 
\end{CD}
$$
is commutative.
Denote by $\Phi_\circ$
the restriction of $\Phi$ to 
the $C^*$-subalgebra $\B_\circ$
of $\OH$ generated by 
$z \in \B_1, w \in \B_2$.
By assumption,
both the maps $\Phi_i: \B_i \longrightarrow  \widehat{O}_\H, i=1,2$
are injective,
so that $\Phi_\circ:\B_\circ \longrightarrow \widehat{O}_\H$
is injective by Lemma \ref{lem:phi0}.
By the preceding lemma,
$
\Phi|_{{\cal F}_\H}: {\cal F}_\H \longrightarrow  \widehat{F}_\H
$
is an isomorphism.
Since the conditional expectation 
${\cal E}_\H: \OH \longrightarrow {\cal F}_\H$
is faithful,
a routine argument shows that $\Phi$ is injective and hence isomorphic.
\end{pf}
Therefore we have
\begin{thm}\label{thm:main2}
For a $C^*$-quad module $\H$
of finite type,
the  $C^*$-algebra 
$\OFH$ 
generated by the quotients 
$[s_\xi], [t_\xi]$
of the creation operators 
$s_\xi, t_\xi$ for $\xi \in \H$
on the Fock spaces $F(\H)$
is canonically isomorphic to the universal $C^*$-algebra
$\OH$
generated by 
operators
$S_1,\dots, S_M, T_1, \dots, T_N$
and elements
$z\in \B_1, w\in \B_2 $
subject to the relations: 
\begin{align*}
\sum_{i=1}^M S_i S_i^* +
\sum_{k=1}^N T_k T_k^* &=1,
\qquad  
S_j^* T_l =0,\\ 
S_i^* S_j  = \langle u_i \mid u_j \rangle_{\B_1}, 
& \qquad
T_k^* T_l  = \langle v_k \mid v_l \rangle_{\B_2}, \\
z S_j = \sum_{i=1}^M S_i  \langle u_i \mid \phi_1(z) u_j \rangle_{\B_1}, 
& \qquad
z T_l = \sum_{k=1}^N T_k  \langle v_k \mid \phi_1(z) v_l \rangle_{\B_2}, \\
w S_j = \sum_{i=1}^M S_i  \langle u_i \mid \phi_2(w) u_j \rangle_{\B_1},
& \qquad
w T_l = \sum_{k=1}^N T_k  \langle v_k \mid \phi_2(w) v_l \rangle_{\B_2},
\end{align*}
for
$
i,j=1,\dots, M, \,
k,l=1,\dots,N$
and
$z \in \B_1, \, w \in \B_2$.
\end{thm}
\begin{pf}  
Theorem \ref{thm:relations} implies that the operators 
$\S_1,\dots, \S_M, \T_1, \dots, \T_N$
and the elements
$\varPhi_1(z), \varPhi_2(w)$ for $z \in \B_1, w\in \B_2$
in $\OFH$ satisfy the eight relations of 
Theorem \ref{thm:main1}.
By Theorem \ref{thm:main1}, we see that 
the correspondences
\begin{equation*}
S_i \longrightarrow \S_i, \quad
T_k \longrightarrow \T_k, \qquad
z \in \B_1 \longrightarrow \varPhi_1(z), \quad
w \in \B_2 \longrightarrow \varPhi_2(w)
\end{equation*}
for 
$i=1,\dots,M,\, k=1,\dots,N,$
and
$z \in \B_1, w\in \B_2$
give rise to an isomorphism from 
$\OH$ to $\OFH$.
\end{pf}
The eight relations of the operators above are called the relations $(\H)$.  
The above generating operators 
$S_1,\dots, S_M$ and $T_1, \dots, T_N$
of the universal $C^*$-algebra $\OH$  
correspond to
two finite basis
$\{u_1,\dots,u_M\}$ and
$\{v_1,\dots,v_N\}$ 
of the Hilbert $C^*$-quad module $\H$
respectively.
On the other hand,
the other $C^*$-algebra
$\OFH$ 
is generated  by the quotients of the creation operators 
$s_\xi, t_\xi$ for $\xi \in \H$
on the Fock spaces $F(\H)$, which 
do not depend on the choice of the two finite bases.
Hence we have
\begin{cor}
\label{cor:cor2}
For a $C^*$-quad module $\H$
of finite type,
the universal $C^*$-algebra
$\OH$
generated by 
operators
$S_1,\dots, S_M, T_1, \dots, T_N$
and elements
$z\in \B_1, w\in \B_2 $
subject to the relations
$(\H)$
does not depend on the choice of the finite bases
$\{u_1,\dots,u_M\}$ and
$\{v_1,\dots,v_N\}$.
\end{cor}

\section{K-Theory formulae}
Let $\H$ be a Hilbert $C^*$-quad module over $(\A;\B_1,\B_2)$
of finite type as in the preceding section.
In this section, we will state K-theory formulae for the $C^*$-algebra $\OFH$.
By the previous section,
the $C^*$-algebra $\OFH$ is regarded as the universal 
$C^*$-algebra $\OH$ generated by 
the operators $S_1,\dots,S_M$ and $T_1,\dots,T_N$
and the elements $z \in \B_1$ and $w \in \B_2$
subject to the relations $(\H)$.
Let us denote by $\B_\circ$ the $C^*$-subalgebra of $\OH$ 
generated by elements $z \in \B_1$ and $w \in \B_2$.
By Lemma \ref{lem:pi}
the correspondence
\begin{equation}
z,w \in \B_\circ \longrightarrow 
\phi_1(z), \phi_2(w) \in C^*(\phi_1(\B_1),\phi_2(\B_2))
\subset {\cal L}_\A(\H)
\end{equation}
gives rise to a $*$-isomorphism
from
$\B_\circ$ onto $C^*(\phi_1(\B_1),\phi_2(\B_2))$
as $C^*$-algebras, which is denoted by 
$\phi_\circ$.
We will restrict our interest to the case when
\begin{enumerate}
\renewcommand{\labelenumi}{(\roman{enumi})}
\item
$S_1, \dots, S_M$ and $T_1, \dots,T_N$ are partial isometries, and 
\item
$S_1 S_1^*, \dots, S_M S_M^*, 
 T_1 T_1^*, \dots, T_N T_N^* $
commute with all elements of 
$\B_\circ$.
\end{enumerate}
If the bases 
$\{u_1,\dots,u_M\}$
and
$\{v_1,\dots,v_N\}$
satisfy the conditions
\begin{equation*}
\langle u_i \mid u_j \rangle_{\B_1} =0 \text{ for } i\ne j,
\qquad
\langle v_k \mid v_l \rangle_{\B_2} =0 \text{ for } k\ne l,
\end{equation*}
the condition (i) holds.
Furthermore if $\phi_1(z)$ acts diagonally on 
$\{u_1,\dots,u_M\}$ for $z \in \B_1$
and
$\phi_2(w)$ acts diagonally on
$\{v_1,\dots,v_N\}$ for $w \in \B_2$,
th condition (ii) holds.
Rcall that the gauge action is denoted by $h$
which is an action of ${\Bbb T}$ 
on $\OH$ such that the fixed point algebra 
$(\OH)^h$ under $h$
is canonically isomorphic to the $C^*$-algebra
$\FH$.
Denote by 
$\widehat{h}$ the dual action of $h$ 
which is an action of ${\Bbb Z} = \widehat{\Bbb T}$
on the $C^*$-crossed product 
$\OH \times_h{\Bbb T}$
by the gauge action $h$ of ${\Bbb T}$.
As in the argument of 
\cite{PasJOT},
$\OH \times_h{\Bbb T}$
is stably isomorphic to 
$\FH$.
Hence we have
$K_*(\OH \times_h{\Bbb T})$
is isomorphic to
$K_*(\FH)$.
The dual action
$\widehat{h}$ induces an automorphism
on the group
$K_*(\OH \times_h{\Bbb T})$
and hence on 
$K_*(\FH)$,
which is denoted by
$\sigma_*$.
Then by \cite{PasJOT} (cf. \cite{Cu2},  \cite{Pim}, etc.)  
we have
\begin{prop}
The following six term
exact sequence of K-theory hold:
\begin{equation*}
\begin{CD}
 K_0(\FH) @>id - \sigma_* >> 
 K_0(\FH) @>\iota_{*}>> 
 K_0(\OH) \\
@AAA @. @VVV \\
 K_1(\OH)
@<<\iota_{*}<
 K_1(\FH) @<<id - \sigma_* < 
 K_1(\FH).  
\end{CD}
\end{equation*}
\end{prop}
We put for $x \in \B_\circ$
\begin{align*}
\lambda_{1,i}(x) & = S_i^* x S_i, \qquad  i=1,\dots,M,\\
\lambda_{2,k}(x) & = T_k^* x T_k, \qquad   k=1,\dots,N.
\end{align*}
Both the families 
$\lambda_{1,i}, \lambda_{2,k}$
yield endomorphisms on $\B_\circ$ 
which give rise to endomorphisms
on the K-groups:
\begin{align*}
\lambda_{1,i *} : &  K_0(\B_\circ) \longrightarrow K_0(\B_\circ),\qquad  i=1,\dots,M,\\\lambda_{2,k *} : &  K_0(\B_\circ) \longrightarrow K_0(\B_\circ),\qquad   k=1,\dots,N.
\end{align*}
We put
$
\lambda_\circ = \sum_{i=1}^M \lambda_{1,i*} + \sum_{k=1}^N \lambda_{2,k*}
$
which is an endomorphism on $K_0(\B_\circ)$.
Now we further assume that
$K_1(\FH) = \{0\}$.
It is routine to show that
the groups
$\Coker(\id - \sigma_*)$ in $K_0(\FH)$
and
$\Ker(\id - \sigma_*)$ in $K_0(\FH)$
are isomorphic to 
the groups
$\Coker(\id - \lambda_\circ)$ in $K_0(\B_\circ)$
and
$\Ker(\id - \lambda_\circ)$ in $K_0(\B_\circ)$
respectively 
by an argument of \cite{Cu2}.
Therefore we have
\begin{prop}
\begin{align*}
K_0(\OH) 
& = \Coker ( \id - \lambda_\circ) \quad \text{ in }\quad  K_0(\B_\circ), \\
K_1(\OH) 
& = \Ker ( \id - \lambda_\circ) \quad \text{ in }\quad  K_0(\B_\circ).
\end{align*}
\end{prop}

\section{Examples}
In this section, we will study the $C^*$-algebras
$\OH$ for the Hilbert $C^*$-quad modules presented in Examples in Section 2.

{\bf 1.}
Let 
$\alpha$,
$\beta$
be automorphisms of a unital 
$C^*$-algebra $\A$
satisfying 
$
\alpha \circ \beta
=
\beta \circ \alpha.
$
Let
${\cal H}_{\alpha, \beta}$
be the associated Hilbert 
$C^*$-quad module of finite type
as in {\bf 1} in Section 2.
It is easy to see 
the following proposition.
\begin{prop} 
The $C^*$-algebra 
${\cal O}_{{\cal H}_{\alpha, \beta}}$
associated to the Hilbert $C^*$-quad module 
${\cal H}_{\alpha, \beta}$
coming from commuting automorphisms
$\alpha,\beta$ of a unital $C^*$-algebra $\A$ 
is isomorphic to the universal $C^*$-algebra 
generated by
two isometries $U,V$ and elements $x$ of $\A$
subject to the following relations:
\begin{gather*}
UU^* + VV^* =1, \\
U U^* x = x U U^*, \qquad 
V V^* x = x V V^*,\\
\alpha(x) = U^* x U, \qquad
\beta(x) = V^* x V
\end{gather*}
for $x \in \A$.
\end{prop}

\medskip

{\bf 2.}
We fix natural numbers $1 < N, M \in {\Bbb N}$.
Consider finite dimensional commutative $C^*$-algebras
$\A = {\Bbb C}, \ \B_1 = {\Bbb C}^N, \ \B_2 = {\Bbb C}^M.$
The algebras $\B_1, \B_2$ have the ordinary product structure and the inner product structure which we denote by 
$\langle \cdot \mid \cdot \rangle_N$ 
and
$\langle \cdot \mid \cdot \rangle_M$ 
respectively.
 Let us denote by
$\H_{M,N}$ the Hilbert $C^*$-quad module ${\Bbb C}^M \otimes {\Bbb C}^N$
over $({\Bbb C};{\Bbb C}^N, {\Bbb C}^M)$
dsefined in {\bf 2} in  Section 2. 
Put the finite bases
\begin{align*}
u_i & = e_i \otimes 1 \in \H_{M,N}, \quad i=1,\dots,M
\text{ as a right } \B_1-\text{module, and}\\
v_k & = 1 \otimes f_k \in \H_{M,N}, \quad k=1,\dots,N
\text{ as a right } \B_2-\text{module.}
\end{align*}
We set
$\Sigma^\circ 
= \{ (i,k) \mid 1 \le i \le M, 1 \le k \le N  \}$
and
put
$e_{(i,k)} =e_i\otimes f_k,\, 
(i,k) \in \Sigma^\circ$
the standard basis of $\H_{M,N}$.
Then the $C^*$-algebra $\B_\circ$
on $\H_{M,N}$ 
generated by
$\B_1$ and $\B_2$ is regarded as 
${\Bbb C}^M \otimes {\Bbb C}^N = \B_2 \otimes \B_1$. 
Hence
$$
\B_\circ = \sum_{(i,k)\in \Sigma^\circ} {\Bbb C}e_{(i,k)}.
$$
\begin{lem}
The $C^*$-algebra 
${\cal O}_{\H_{M,N}}$ is generated by 
operators
$S_i,T_k, e_{(i,k)},i=1,\dots,M, k=1,\dots,N$
satisfying
\begin{align}
\sum_{i=1}^M S_i S_i^*  + & \sum_{k=1}^N T_k T_k^* =1, \\
S_i^* S_j = \delta_{i,j}, \qquad &  \qquad T_k^* T_l = \delta_{k,l}, \\
e_{(i,k)} S_j  = \delta_{i,j} \sum_{h=1}^M S_j e_{(h,k)},
& \qquad
e_{(i,k)} T_l  = \delta_{k,l} \sum_{m=1}^N T_l e_{(i,m)} \label{eqn:eT}   
\end{align}
for $i,j =1,\dots,M,\, k,l=1,\dots,N.$
\end{lem}
\begin{pf}
It suffices to show the equalities \eqref{eqn:eT}.
We have
\begin{align*}
e_{(i,k)} S_j 
& = S_j \langle u_j \mid \phi_\circ(e_{(i,k)}) u_j \rangle_{\B_1} \\
& = S_j \langle e_j \otimes 1 
\mid (e_i \otimes f_k )(e_j\otimes 1) \rangle_{\B_1} \\
& = \delta_{i,j}  S_j (1 \otimes f_k)  
  = \delta_{i,j}  S_j \sum_{h=1}^M e_{(h,k)}.
\end{align*}
The other equality of
\eqref{eqn:eT} is similarly shown.
\end{pf}
Put
$$
S_{(i,k)} = e_{(i,k)} S_i, \qquad
T_{(i,k)} = e_{(i,k)} T_k
\qquad \text{ for } (i,k) \in \Sigma^\circ.
$$
Then we have
\begin{lem}
\begin{align}
e_{(i,k)} & = S_{(i,k)}S_{(i,k)}^* + T_{(i,k)} T_{(i,k)}^*, \label{eqn:(a)}\\
 S_i = \sum_{k=1}^N S_{(i,k)}, & \qquad
 T_k = \sum_{i=1}^M  T_{(i,k)}, \label{eqn:(b)} \\
\sum_{(i,k) \in \Sigma^\circ}S_{(i,k)}S_{(i,k)}^* 
& + \sum_{(i,k) \in \Sigma^\circ} T_{(i,k)} T_{(i,k)}^* =1, \label{eqn:(c)}\\
S_{(i,k)}^*S_{(i,k)} 
&= \sum_{j=1}^M (S_{(j,k)}S_{(j,k)}^* +T_{(j,k)} T_{(j,k)}^* ),
\label{eqn:(d)S} \\
T_{(i,k)}^*T_{(i,k)} 
&= \sum_{l=1}^N (S_{(i,l)}S_{(i,l)}^* +T_{(i,l)} T_{(i,l)}^* )  
\label{eqn:(d)T} 
\end{align}
for $i=1,\dots, M, \, k=1,\dots,N$ and
$(i,k) \in \Sigma^\circ.$
\end{lem}
\begin{pf}
Since
$e_{(i,k)} S_j = \delta_{i,j}e_{(i,k)} S_i$,
we have
$$ 
S_{(i,k)}S_{(i,k)}^* 
= e_{(i,k)} S_i S_i^* e_{(i,k)} 
= e_{(i,k)}(\sum_{j=1}^M S_j S_j^*) e_{(i,k)} 
$$
and similarly
$ 
T_{(i,k)}T_{(i,k)}^* 
= e_{(i,k)}(\sum_{l=1}^N T_l T_l^*) e_{(i,k)}.
$
Hence we have
$$
e_{(i,k)}  
= 
e_{(i,k)}(\sum_{j=1}^M S_j S_j^* 
+ \sum_{l=1}^N T_l T_l^*) e_{(i,k)}
=S_{(i,k)}S_{(i,k)}^* + T_{(i,k)} T_{(i,k)}^*,
$$
so that \eqref{eqn:(a)} holds.
As
$1 = \sum_{(j,k) \in \Sigma^\circ} e_{(j,k)}$,
the equality \eqref{eqn:(c)} holds.
Since
$e_{(j,k)} S_i = 0$ for  $j \ne i$,
we have
$$
S_i = (\sum_{(j,k) \in \Sigma^\circ} e_{(j,k)})S_i
 = \sum_{k=1}^N S_{(i,k)}
$$
and similarly 
$ T_k = \sum_{i=1}^M  T_{(i,k)}$,
so that \eqref{eqn:(b)} holds.
By \eqref{eqn:eT},
it follows that
\begin{equation*}
S_{(i,k)}^*S_{(i,k)} 
= S_i^* \sum_{j=1}^M S_i e_{(j,k)} 
= \sum_{j=1}^M e_{(j,k)} 
= \sum_{j=1}^M (S_{(j,k)}S_{(j,k)}^* +T_{(j,k)} T_{(j,k)}^* ),
\end{equation*}
and smilarly we have
$$
T_{(i,k)}^*T_{(i,k)} 
= \sum_{l=1}^N (S_{(i,l)}S_{(i,l)}^* +T_{(i,l)} T_{(i,l)}^* ).  
$$
\end{pf}
\begin{thm}\label{thm:CK} 
The $C^*$-algebra ${\cal O}_{\H_{M,N}}$
associated with the 
Hilbert $C^*$-quad module 
$\H_{M,N} ={\Bbb C}^M \otimes {\Bbb C}^N$
is 
generated by partial isometries
$
S_{(i,k)}, T_{(i,k)}
$ 
for
$(i,k)\in \Sigma^\circ =\{(i,k)\mid i=1,\dots,M, k=1,\dots,N \}$
satisfying the  relations: 
\begin{align*}
\sum_{(i,k) \in\Sigma^\circ} 
& S_{(i,k)}S_{(i,k)}^*
+
\sum_{(i,k) \in\Sigma^\circ} 
T_{(i,k)}T_{(i,k)}^* = 1, \\
S_{(i,k)}^*S_{(i,k)}
=
&
\sum_{j=1}^M
S_{(j,k)}S_{(j,k)}^* + T_{(j,k)} T_{(j,k}^*),\\
T_{(i,k)}^* T_{(i,k)}
=
&
\sum_{l=1}^N 
(
S_{(i,l)}S_{(i,l)}^* + T_{(i,l)} T_{(i,l)}^*)
\end{align*}
for $(i,k) \in  \Sigma^\circ$.
\end{thm}
\begin{pf}
By the preceding lemma, one knows that 
$e_{(i,k)}, S_i, T_k$ are generated by 
the operators
$S_{(i,k)}, T_{(i,k)}$
so that the algebra
${\cal O}_{\H_{M,N}}$ is generated by 
the partial isometries $S_{(i,k)}, T_{(i,k)},
(i,k) \in\Sigma^\circ$.
\end{pf}
Let
$I_n$ be the $n \times n$ identity matrix
and 
$E_n$ the $n \times n$ matrix whose entries are all $1'$s.
For an $M \times M$-matrix
$C = [c_{i,j}]_{i,j=1}^M$
and
an $N \times N$-matrix
$ D = [d_{k,l}]_{k,l=1}^N$,
denote by
$C \otimes D$
 the $MN \times MN$  matrix
 $$
C \otimes D
=
\begin{bmatrix}
c_{11}D & c_{12}D  & \dots & c_{1M}D  \\
c_{21}D & c_{22}D  & \dots & c_{2M}D  \\
\vdots  & \vdots   & \ddots& \vdots   \\
c_{M1}D & c_{M2}D  & \dots & c_{MM}D  
\end{bmatrix}
$$
so that
\begin{equation*}
E_M \otimes I_N =
\begin{bmatrix}
I_N    & I_N    & \dots & I_N    \\
I_N    & I_N    & \dots & I_N    \\
\vdots & \vdots & \ddots& \vdots \\
I_N    & I_N    & \dots & I_N 
\end{bmatrix},
\qquad
I_M \otimes E_N =
\begin{bmatrix}
E_N   & 0    & \dots & 0      \\
0     & E_N  & \ddots & \vdots \\
\vdots &\ddots      & \ddots& 0      \\
0     & \dots& 0     & E_N 
\end{bmatrix}.
\end{equation*}
The index set  
$\{(i,k) \mid i=1,\dots,M, k=1,\dots,N \}$
of the standard basis of ${\Bbb C}^M \otimes {\Bbb C}^N$ 
is ordered  lexicographically  from left as in the following way:
\begin{equation*}
 (1,1),\dots, (1, N),
   (2,1),\dots, (2,N),
\dots,
(M,1),\dots, (M,N). 
\end{equation*}
Put the $MN \times MN$ matrices
\begin{equation*}
A_{M,N} = E_M \otimes I_N, \qquad
B_{M,N} = I_M \otimes E_N
\end{equation*}
and the $2MN \times 2MN$ matrix
\begin{equation*}
H_{M,N} =
\begin{bmatrix}
A_{M,N} & A_{M,N} \\
B_{M,N} & B_{M,N}   
\end{bmatrix}.
\end{equation*}
Then we have
\begin{thm}
The $C^*$-algebra 
${\cal O}_{\H_{M,N}}$
is isomorphic to the Cuntz-Krieger algebra
${\cal O}_{H_{M,N}}$ for the matrix  
$H_{M,N}$.
The algebra ${\cal O}_{H_{M,N}}$
is simple and purely infinite and
is  isomorphic to the Cuntz-Krieger algebra 
${\cal O}_{A_{M,N} + B_{M,N}}$
for the matrix 
$A_{M,N} + B_{M,N}$.
\end{thm}
\begin{pf}
By the preceding proposition,
the $C^*$-algebra 
${\cal O}_{\H_{M,N}}$
is isomorphic to the Cuntz-Krieger algebra
${\cal O}_{H_{M,N}}$ for the matrix  
$H_{M,N}$.
Since the matrix 
$H_{M,N}$ is aperiodic,
the algebra is simple and purely infinite.
The $n$-th column of  
the matrix 
$H_{M,N}$ 
coincides with the
$n+N$-th column for every 
$n=1,\dots,M$.
One sees that
the matrix $A_{M,N} + B_{M,N}$
is obtained from $H_{M,N}$
by amalgamating them.
The precedure is called the column amalgamation
and induces an isomorphism on their Cuntz-Krieger algebras
(see \cite{MaPre2012}). 
\end{pf}
In \cite{MaPre2012},
the abelian groups
$$
{\Bbb Z}^{MN}/ (A_{M,N} + B_{M,N} - I_{MN}){\Bbb Z}^{MN},
\qquad
\Ker(A_{M,N} + B_{M,N} - I_{MN}) \text{ in }{\Bbb Z}^{MN}
$$
have been computed by using Euclidean algorithms.
For the case 
$M=2$, they are
$$
{\Bbb Z}/(N^2 -1){\Bbb Z}, \qquad
\{0\}
$$
respectively,
so that
we see 
$$
K_0({\cal O}_{\H_{2,N}}) = {\Bbb Z}/(N^2 -1){\Bbb Z},
\qquad
K_1({\cal O}_{\H_{2,N}}) = 0
$$
(see \cite{MaPre2012} for details).

\medskip

{\bf 3.}
For a $C^*$-textile dynamical system
$(\A,\rho,\eta, \Sigma^\rho,\Sigma^\eta,\kappa)$,
let $\H^{\rho,\eta}_\kappa$ 
be the $C^*$-quad module over $(\A;\B_1,\B_2)$
as in {\bf 3} in Section 2.
The $C^*$-algebra
${\cal O}_{\H^{\rho,\eta}_\kappa}$
has been studied in \cite{MaPre20112}.

\section{Higher dimensional analogue} 
In this final section,
we will state a generalization of Hilbert $C^*$-quad modules
to Hilbert modules with multi actions of $C^*$-algebras. 

Let $\A$ be a unital $C^*$-algebra
and 
$\B_1,\dots,\B_n$ be $n$-family of unital $C^*$-algebras.
Suppose that there exists a unital embedding
\begin{equation*}
\iota_i: \A \hookrightarrow \B_i
\end{equation*}
for each $i=1,\dots,n$.
Suppose that there exists
 a right action $\psi_i$
of $\A$ on $\B_i$ such that
\begin{equation*}
b_i \psi_i(a) \in \B_i \qquad
\text{ for } \quad 
b_i \in \B_i, \, a \in \A, \, i=1,\dots,n. 
\end{equation*} 
Hence $\B_i$ is a right $\A$-module through $\psi_i$ for $i=1,\dots,n$.
Let 
$\H$ be a Hilbert $C^*$-bimodule over $\A$
with a right action of $\A$, an $\A$-valued 
right inner product $\langle \cdot \mid \cdot \rangle_\A$
and a $*$-homomorphism 
$\phi_\A$ from $\A$ to ${\cal L}_\A(\H)$.
It is called 
 a Hilbert $C^*$-{\it multi module over} 
$(\A;\B_i,i=1,\dots,n)$
if  $\H$ has a multi structure of Hilbert $C^*$-bimodules over $\B_i$
for $i=1,\dots,n$
such that for each $i=1,\dots,n$ there exist
a right action $\varphi_i$ of $\B_i$ on $\H$ 
and 
a left action $\phi_i$ of $\B_i$ on $\H$
and 
a $\B_i$-valued right inner product 
$
\langle \cdot \mid \cdot \rangle_{\B_i}
$
such that 
$ 
\phi_i(z_i)  \in {\cal L}_\A(\H)
$ 
and
\begin{equation*}
[\phi_i(z_i) \xi]\varphi_j(w_j)  = \phi_i(z_i)[ \xi \varphi_j(w_j) ],  
\qquad
 \xi \varphi_j(z_j \psi_j(a))  = [\xi \varphi_j(z_j)]a 
\end{equation*}
for  $\xi\in \H, z_i\in \B_i, w_j \in \B_j, a \in \A$,
$i,j=1,\dots,n$
and
\begin{equation*}
\phi_\A(a) = \phi_i(\iota_i(a)), \qquad a \in \A, \, i=1,\dots,n. 
\end{equation*}
The operator 
$\phi_i(z_i)$ on $\H$ is adjointable
with respect to the inner product
$\langle \cdot \mid \cdot \rangle_{\B_i}$
whose adjoint 
$\phi_i(z_i)^*$ 
coincides with the adjoint of $\phi_i(z_i)$
with respect to the inner product
$\langle \cdot \mid \cdot \rangle_\A$
so that
$\phi_i(z_i)^* = \phi_i(z_i^*)$. 
We assume that the left actions $\phi_i$
of $\B_i$ on $\H$ for $i=1,2$ are faithful.
We require the following compatibility conditions between 
the right $\A$-module structure of $\H$ and 
the right $\A$-module structure of $\B_i$
through $\psi_i$:
\begin{equation*}
\langle \xi \mid \eta a \rangle_{\B_i}
=\langle \xi \mid \eta  \rangle_{\B_i} \psi_i(a),
\qquad
\xi, \eta \in \H, a \in \A, \, i=1,\dots,n. 
\end{equation*} 
We further 
assume that $\H$ is a full Hilbert $C^*$-bimodule with respect to 
the inner product
$
\langle \cdot \mid \cdot \rangle_\A,
\langle \cdot \mid \cdot \rangle_{\B_i}
$
for each.
A Hilbert $C^*$-multi module $\H$ over $(\A;\B_i,i=1,\dots,n)$
is said to be {\it of general type}
if there exists a faithful completely positive map
$\lambda_i: \B_i \longrightarrow \A$ for $ i=1,\dots,n$ 
such that
\begin{align*}
\lambda_i(b_i \psi_i(a))&  = \lambda_i(b_i) a,
\qquad  b_i \in \B_i, a \in \A,  \\ 
\lambda_i(\langle \xi \mid \eta \rangle_{\B_i}) 
& = 
\langle \xi \mid \eta \rangle_{\A}, \qquad
\xi, \eta \in \H, \, i=1,\dots,n.
\end{align*}
A Hilbert $C^*$-multi module $\H$ over $(\A;\B_i,i=1,\dots,n)$
is said to be {\it of finite type}
if there exists a family  
$\{ u^{(i)}_1,\dots,u^{(i)}_{M^{(i)}} \}, i=1,\dots,n$ 
of finite bases
of $\H$ as a right Hilbert $\B_i$-module
for each $i=1,\dots,n$ 
such that
\begin{align*}
\sum_{j=1}^{M^{(i)}} & u^{(i)}_j \varphi_i(\langle u^{(i)}_j \mid \xi \rangle_{\B_i})
= \xi,
\qquad \xi \in \H, \, i=1,\dots,n \\
\intertext{and}
\langle u^{(i)}_j \mid  & \phi_k(w_k) u^{(i)}_h \rangle_{\B_i} 
 \in \A, \qquad   w_k \in \B_k, \, j,h = 1,\dots,M^{(i)},  \\ 
\sum_{j=1}^{M^{(i)}} & \langle u^{(i)}_j \mid 
\phi_k(\langle \xi \mid \eta \rangle_{\B_k}) u^{(i)}_j \rangle_{\B_i}
 = \langle \xi \mid \eta \rangle_\A 
\end{align*}
for all $\xi,\eta \in \H$, $i,k=1,\dots,n$
with $i \ne k$.

By a generalizing argument to the preceding sections,
we may construct a $C^*$-algebra $\OFH$ associated with
the Hilbert $C^*$-multi module $\H$ 
by a similar manner to the preceding sections,
that is, the $C^*$-algebra generated by $n$-kinds of creation operators
$s^{(i)}_\xi, \xi \in \H, i=1,\dots,n$ 
on the generalized Fock space $F(\H)$
by the ideal generated by the finite rank operators.
One may show the following generalization:

\begin{prop}
Let $\H$ be a Hilbert $C^*$-multi module over
$(\A;\B_i,i=1,\dots,n)$
of finite type
with
a finite basis $\{ u^{(i)}_1,\dots,u^{(i)}_{M^{(i)}} \}$ of $\H$ 
as a Hilbert $C^*$-right module over $\B_i$
for each $i=1,\dots,n$.
Then
the  $C^*$-algebra 
$\OFH$ 
generated by the $n$-kinds of creation operators 
on the generalized Fock spaces $F(\H)$
is canonically isomorphic to the universal $C^*$-algebra
$\OH$
generated by the operators
$S^{(i)}_1,\dots, S^{(i)}_{M^{(i)}}$ 
and elements
$z_i\in \B_i$
for $i=1,\dots,n$
subject to the relations: 
\begin{gather*}
\sum_{i=1}^n  \sum_{k=1}^{M^{(i)}} S^{(i)}_k S^{(i)*}_k  =1, 
 \qquad
S^{(i)*}_k S^{(j)}_m =0, \quad i\ne j, \\
S^{(i)*}_k S^{(i)}_l = \langle u^{(i)}_k \mid u^{(i)}_l \rangle_{\B_i}, 
\qquad
z_j S^{(i)}_k 
 = \sum_{l=1}^{M^{(i)}} S^{(i)}_l \langle u^{(i)}_l \mid \phi_j(z_j) u^{(i)}_k \rangle_{\B_i} 
\end{gather*}
for
$
z_j \in \B_j, \, i,j=1,\dots, n,\,
k,l=1,\dots,M^{(i)}, \, 
m=1,\dots,M^{(j)}.
$
\end{prop}
The proof of the above proposition is similar to the proof of Theorem 1.1.


\end{document}